\documentclass[12pt,a4paper]{article}

\usepackage[english,francais]{babel} 
\RequirePackage[latin1]{inputenc}
\RequirePackage[T1]{fontenc}
\RequirePackage{graphicx,hyperref,ifpdf,color}

\RequirePackage{amsmath,amssymb,amsfonts,bbm}

\newcommand{\bs}{{\bf s}}

\newcommand{\bT}{{\bf T}}

\newcommand{\R}{\mathbb{R}}
\newcommand{\N}{\mathbb{N}}

\renewcommand{\P}{\mathbb{P}}

\newcommand{\E}{\mathbb{E}}
\newcommand{\C}{\mathbb{C}}
\newcommand{\D}{\mathbb{D}}

\newcommand{\TT} { {\cal T }}

\newcommand{\FF} { {\cal F }}

\def\build#1_#2^#3{\mathrel{
\mathop{\kern 0pt#1}\limits_{#2}^{#3}}}

\def\cq{$\hfill \square$}
\def\un{\underline}
\def\d{{\rm d}}

\def\eps{\varepsilon}

\def\ba{\begin{eqnarray*}}
\def\ea{\end{eqnarray*}}
\def\ov{\overline}

\def\wt{\widetilde}

\newcommand{\ind}{\mathbbm{1}}

\newtheorem{thm}{Theorem}
\newtheorem{lmm}{Lemma}
\newtheorem{prp}{Proposition}
\newtheorem{defn}{Definition}

\usepackage{hyperref} 
\usepackage{verbatim,multicol}

\def\miermont{{\href{http://www.dma.ens.fr/~miermont/}{Gr\'egory Miermont}}}
\def\dma{{\href{http://www.dma.ens.fr/}{DMA}}}
\def\lpma{{\href{http://www.proba.jussieu.fr/}{LPMA}}}

\begin{document}
\title{Self-similar fragmentations derived 
from the stable tree I: 
splitting at heights}
\author{\miermont
\\ \dma, \'Ecole Normale Sup\'erieure, \\
and \lpma, Universit\'e Paris VI.\\
45, rue d'Ulm, \\
75230 Paris Cedex 05\thanks{Research supported in part by NSF Grant
DMS-0071448}}
\date{}

\maketitle

\begin{abstract}
The basic object we consider is a certain model of continuum random tree, 
called the stable tree. We construct a fragmentation process 
$(F^-(t),t\geq 0)$ out of this
tree by removing the vertices located under height $t$. 
Thanks to a self-similarity property of the stable tree, we show that 
the fragmentation process is also self-similar. The semigroup and 
other features of the fragmentation are given explicitly. Asymptotic results 
are given, as well as a couple of related results on continuous-state 
branching processes. 
%As proved in a companion paper, another method for fragmenting the 
%stable tree induces another self-similar fragmentation with same 
%characteristics as the ones considered here, except for the speed at which 
%fragments decay. 
\end{abstract}


\bigskip

\noindent{\bf A.M.S. Classification. }60J25, 60G52, 60J80.
\newpage

\section{Introduction}

The recent advances in the study of coalescence and fragmentation processes
pointed at the key role played by tree structures in this topic, both
at the discrete and continuous level \cite{jpse96cmc,jpda98sac,jpda97ebac}. 
Our goal here is to push further the investigation, begun in
\cite{jpda98sac,bertsfrag02}, of a category of fragmentations
obtained by cutting a certain class a continuum random tree. The tree that was
fragmented 
in the latter articles is the Brownian Continuum Random Tree of Aldous,
and the fragmentation is related to the so-called {\em standard additive
coalescent}. The family of trees we consider is a natural but 
technically involved ``L\'evy generalization'' of the
Brownian tree. It has been introduced in Duquesne and Le Gall 
\cite{duqleg02}, and implicitly considered in the previous work of Kersting 
\cite{kersting98m}. Some of these trees, 
as their Brownian companion, enjoy certain
self-similar properties. In the present work the crucial property is that  
when removing the vertices of the stable tree located under a fixed 
{\em height} (or distance to the root), 
the remaining object is a forest of smaller 
trees that have the same law as the original one up to rescaling.
This is formalized in Lemma \ref{spropfm} below. 
This way of logging the stable tree induces a {\em fragmentation process} 
which by the property explained above turns out to be a
{\em self-similar fragmentation}, the theory
of such processes being extensively studied by Bertoin 
\cite{berthfrag01,bertsfrag02,bertafrag02}. 
The goal of this paper is to describe the characteristics and give some
properties of this
fragmentation process. We will have to use stochastic processes and 
combinatorial approaches in the same time; in particular, we will encounter 
$\sigma$-finite generalizations of the $(\alpha,\theta)$-partitions of 
\cite{pitmancsp02}, which are distributions on the set of partitions
of $\N=\{1,2,\ldots\}$, as well as we will need the construction of the 
stable tree out of L\'evy processes and its connection to continuous-state
branching processes (CSBP) explained in \cite{duqleg02}. 

In a companion paper \cite{mierfplus} we will consider another way of obtaining
a self-similar fragmentation by another cutting device on the stable tree, 
using the heuristic fact that when cutting at random one hub in the 
the stable tree, the trunk and branches that have been separated are scaled 
versions of the initial tree. 
Surprisingly, although this other device looks quite different from the first
(no mass is lost when cutting a hub, whereas there is a loss of mass when 
we throw everything that is located under the height $h$), it turns out that 
the only difference between these two fragmentations is the speed at which 
fragments decay. 

To state our main results, let us introduce quickly the already mentioned
tree structures and fragmentation processes, postponing the details to 
a further section. 

Let $S=\{\bs=(s_1,s_2,\ldots):s_1\geq s_2\geq\ldots\geq 0, \sum_{i\geq 1}
s_i\leq 1\}$. A ranked self-similar
fragmentation process $(F(t),t\geq 0)$ with index $\beta\in\R$ 
is a $S$-valued Markov process that is continuous in probability,
such that $F(0)=(1,0,0,\ldots)$ and such that conditionally on 
$F(t)=(x_1,x_2,\ldots)$, $F(t+t')$ has the law of the decreasing arrangement
of the sequences $x_i F^{(i)}(x_i^{\beta}t')$, where the $F^{(i)}$ are 
independent with the same law as $F$. That is, after time $t$, the 
different fragments evolve independently with a speed that depends on their 
size. It has been shown in \cite{bertsfrag02} that such fragmentations are
characterized by a 3-tuple $(\beta,c,\nu)$, where $\beta$ is 
the index, $c\geq 0$ is an ``erosion'' real constant saying that the fragments
may melt continuously at some rate depending on $c$, and $\nu$ 
is a $\sigma$-finite
measure on $S$ that attributes mass $0$ to $(1,0,\ldots)$ and that 
integrates $\bs \mapsto (1-s_1)$. This measure governs the sudden dislocations
in the fragmentation process, and the integrability assumption ensures that 
these dislocations do not occur too quickly, although the fragmentation
epochs may form a dense subset of $\R_+$ as soon as $\nu(S)=+\infty$.
When $\beta<0$, a positive fraction of the mass can disappear 
within a finite time, even though there is no loss of mass due to 
erosion nor to sudden dislocations. This phenomenon will be 
crucial in the fragmentation $F^-$ below.

The trees we are considering are continuum random trees. Intuitively, they are
metric spaces 
with an ``infinitely ramified'' tree structure, which can be
considered as genealogical structures combined with two measures: a
$\sigma$-finite {\em length measure} supported by the ``skeleton'' of  the tree
and a finite {\em mass measure} supported by its leaves, which are everywhere
dense in the tree. 
These trees can be defined in several equivalent ways:
\begin{itemize}
\item as a weak limit of Galton-Watson trees
\item through its {\em height process} $H$, which is a positive continuous
process on $[0,1]$. To a point $u\in[0,1]$ corresponds a vertex of the 
tree with height (distance to the root) equal to $H_u$, and the mass measure
on the tree is represented by Lebesgue's measure on $[0,1]$
\item through its explicit ``marginals'', that is, the laws of subtrees spanned
by a random sample of leaves.
\end{itemize}
We will have to use the second (stochastic process)
and third (combinatorial) points of view. We know from
the works of Duquesne and Le Gall \cite{duqleg02} and Duquesne \cite{duq02}
that one may define a particular instance of tree, called the 
{\em stable tree with index $\alpha$} (for some $\alpha\in(1,2]$).
When $\alpha=2$, the stable tree is equal to the Brownian 
CRT of Aldous \cite{aldouscrt93}. 
We will recall the rigorous construction of the height 
process of the stable tree in Sect.\ \ref{secsttree},  but let us state our
results now. Fix $\alpha\in(1,2)$ and let 
$(H_s,0\leq s\leq 1)$ be the height process of the stable tree with index 
$\alpha$. 

The fragmentation process, that we call $F^-$, is defined
as follows. 
For each $t\geq0$, let $I_-(t)$ be the open subset of $(0,1)$ defined by
$$I_-(t)=\{s\in(0,1) : H_s>t\}.$$
With our intuitive interpretation of the height process, 
$I_-(t)$ is the set of vertices of the tree with height $>t$.
We denote by $F^-(t)$ the decreasing sequence of the lengths of the connected
components of $I_-(t)$. Hence, $F^-(t)$ is the sequence of the masses of
the tree components obtained by cutting the stable tree below height $t$.
Notice that $F^-$ is a direct generalization of the fragmentation 
$F$ in \cite[Section 4]{bertsfrag02}. The boundedness of $H$ implies that
$F^-(t)=(0,0,\ldots)$  as soon as $t\geq \max_{0\leq s\leq 1} H_s$.

\begin{prp}\label{P1}
The process $F^-$ is a ranked self-similar fragmentation with 
index $1/\alpha-1\in(-1/2,0)$ and erosion coefficient $0$. 
\end{prp}

Notice that, as mentioned before, $F^-$ loses some mass, and eventually 
disappears completely in finite time even though the erosion is $0$. 
This is due, of course, to the fact that the self-similarity index is 
negative. 

Our main result is a description of the dislocation 
measure $\nu_-(\d\bs)$ of $F^-$.  
Let us introduce 
some notation. For $\alpha\in(1,2)$, let $(T_x,x\geq 0)$ 
be a stable subordinator with Laplace exponent $\lambda^{1/\alpha}$, 
that is, $T_x$ is the sum of the magnitudes of the atoms of
a Poisson point process on $(0,\infty)$ with intensity $c_{\alpha}x \d r/
r^{1+1/\alpha}$, where $c_{\alpha}=(\alpha\Gamma(1-1/\alpha))^{-1}$. 
%Let $(q^{(1/\alpha)}_x(t),t\geq 0)$ be the density of $T^{(1/\alpha)}_x$.
We denote by $\Delta T_x=T_x-T_{x-}$ 
the jump at level $x$ and by $\Delta T_{[0,x]}$ 
the sequence of the jumps of $T$ before time $x$, 
and ranked in decreasing order. 
%Results of Perman \cite{perm93} allow
%us to consider such sequences of jumps conditioned on their sum 
%$T^{(1/\alpha)}_x$. We will give it in detail later in section \ref{levyproc},
%for now we state our main result:
Define the measure $\nu_{\alpha}$ on $S$ by 
\begin{equation}\label{nupm}
\nu_{\alpha}(\d\bs)=E\left[T_1\, ;\, \frac{
\Delta T_{[0,1]}}{T_1}
\in \d\bs\right]
%\int_{\R_+}\d \ell \;\frac{q^{(1/\alpha)}_{\ell}(1)}{\ell^{\alpha+1}}
%P\left((\Delta T^{(1/\alpha)}_{[0,\ell]})\in \d\bs\left|
%T^{(1/\alpha)}_{\ell}=1\right.\right)\; ,
\end{equation}
where the last expression means 
that for any positive measurable function $G$, the quantity
$\nu_{\alpha}(G)$ is equal to $E[T_1 \, G(T_1^{-1} 
\Delta T_{[0,1]})]$.

\begin{thm}\label{T1}
The dislocation measure of $F^-$ is $\nu_-=D_{\alpha}\nu_{\alpha}$, where
$$D_{\alpha}=\frac{\alpha(\alpha-1)\Gamma\left(1-\frac{1}{\alpha}\right)
}{\Gamma(2-\alpha)}
=\frac{\alpha^2\Gamma\left(2-\frac{1}{\alpha}\right)}{\Gamma(2-\alpha)}.$$
\end{thm} 

Some comments about this. First, the dislocation measure 
charges only the sequences 
$\bs$ for which $\sum_{i\geq 1}s_i=1$, that is, no mass can be lost within 
a sudden dislocation. 
%Second, the scaling properties of $T^{(1/\alpha)}$, and 
%its density $q_{\ell}$, given by 
%\begin{equation}\label{scalprT}
%(\ell^{-\alpha}T_{x\ell},x\geq 0)\build=_{}^{d}
%(T_x,x\geq 0)\, ,\quad q_{\ell}(1)=\ell^{-\alpha}q_1(\ell^{-
%\alpha}),
%\end{equation}
%allows us to write $\nu^{\alpha}_{\pm}$ in the form
%$$\nu^{\alpha}_{\pm}(\d\bs)=D_\alpha\int_{\R_+}\d \ell \frac{q_1(
%\ell^{-\alpha})}{\ell^{2\alpha-1}}P\left(\ell^{\alpha}
%(\Delta T_{[0,1]})\in \d\bs
%\left|T_1=\ell^{-\alpha}
%\right.\right).$$
%Letting $u=\ell^{-\alpha}$ gives
%\begin{eqnarray}
%\nu^{\alpha}_{\pm}(\d\bs)&=&\frac{D_{\alpha}}{\alpha}
%\int_{\R_+}\d u\; u\, q_1(u)
%P\left(u^{-1}(\Delta T_{[0,1]})\in \d\bs
%\left|T_1=u\right.\right)\nonumber\\
%&=&\frac{D_{\alpha}}{\alpha}E\left[T_1\, ;\, \frac{
%(\Delta T_{[0,1]})}{T_1}
%\in \d\bs\right],\label{nupmt}
%\end{eqnarray}
%with obvious notations. 
Second, we recognize an expression close to \cite{pitmanpk02}, 
of a Poisson-Dirichlet type
distribution. However, it has to be noticed that this corresponds to 
a forbidden parametrization $\theta=-1$, and indeed, the measure that 
we obtain is infinite since $E[T_1]=\infty$. 
This measure integrates $1-s_1$ though, just as it has to. Indeed, 
$E[T_1-\Delta_1]$ is finite if $\Delta_1$ denotes the largest
jump of $T$ before time $1$. To see this, notice that
$\Delta_1\geq \Delta^*_1$ where $\Delta^*_1$ 
is a size-biased pick from the jumps
of $T$ before time $1$, and it follows from Lemma \ref{permanform} in Sect.\
\ref{levyproc} below and scaling arguments that $T-\Delta^*_1$ has finite
expectation. 

The rest of the paper is organized as follows. 
In Sect.\ \ref{prelim} we first recall some facts about L\'evy processes,
excursions, and conditioned subordinators. Then we 
give the rigorous description of the stable tree, and 
state some properties of the height process that 
we will need. Last we recall some facts about 
self-similar fragmentations. We then obtain the characteristics of 
$F^-$ in Sect.\ \ref{fmoins} and derive its semigroup. 
We insist on the fact that knowing explicitly the semigroup of 
a fragmentation process is in general a very complicated problem, 
see \cite{MS03} for somehow surprising negative results in this vein. 
However, most of 
the fragmentation processes that have been extensively studied
in recent years \cite{jpda98sac,bertfrag99,mier01,bertsfrag02} do have known, 
and sometimes strange-looking semigroups involving conditioned Poisson
clouds. And as a matter of fact, the fragmentation $F^+$ considered in the 
companion paper \cite{mierfplus} has also an explicit semigroup. 
We end the study of $F^-$ by 
giving asymptotic results for small times in Sect.\ \ref{asymzero}. These 
results need some properties of conditioned continuous-time branching 
processes, which are in the vein of Jeulin's results for the 
rescaled Brownian excursion and its local times.
We prove these properties in 
Sect.\ \ref{secCSBP}, where we give the rigorous definition of some processes
that are used heuristically in Sect.\ \ref{fmoins} to conjecture the 
form of the dislocation measure. 

\section{Preliminaries}\label{prelim}

\subsection{Stable processes, excursions, conditioned inverse subordinator}
\label{levyproc}

Throughout the paper, 
we let $(X_s,s\geq 0)$ be the canonical process in the Skorokhod space 
$\D([0,\infty))$ of c\`adl\`ag paths on $[0,\infty)$. 
Recall that a L\'evy process is a real-valued c\`adl\`ag process with
independent and stationary increments. We fix
$\alpha\in(1,2)$. Let $P$ be the law that makes $X$ a stable L\'evy process
with no negative jumps and Laplace exponent 
$E[\exp(-\lambda X_s)]=\exp(\lambda^{\alpha})$ for $s,\lambda\geq 0$,
where $E$ is the expectation associated with $P$.  Such a process has infinite
variation and satisfies $E[X_1]=0$. When there is no ambiguity, we may
sometimes speak of $X$ as being itself the  L\'evy process with law $P$.
Writing this in the form of the L\'evy-Khintchine formula, we have :
\begin{equation}\label{laplaX}
E[\exp(-\lambda
X_s)]=\exp\left(s\int_0^{\infty}\frac{C_{\alpha}\d x}{x^{1+\alpha}}(e^{-\lambda
x}-1+\lambda x)\right),\quad s,\lambda\geq 0,
\end{equation} 
where $C_{\alpha}=\alpha(\alpha-1)/\Gamma(2-\alpha)$. In particular, the L\'evy
measure of $X$ under $P$ is $C_{\alpha}x^{-1-\alpha}\d x \ind_{\{x>0\}}$.
%Let $(\FF_s,s\geq 0)$ be the usual
%completed natural filtration of $X$ under $P$. 
An important property of $X$ is then the scaling property: under $P$,
$$\left(\frac{1}{\lambda^{1/\alpha}}X_{\lambda s},s\geq 0\right)
\build=_{}^{d}(X_{s},s\geq 0) \quad \mbox{ for all } \lambda>0.$$
It is known \cite{skorokhod} that under $P$, $X_s$ has a density 
$(p_s(x),x\in \R)$ for every $s>0$, such that $p_s(x)$ is jointly 
continuous in $x$ and $s$.

\paragraph{Excursions}
Let $\un{X}$ be the infimum process of $X$, defined for $s\geq 0$ by
$$\un{X}_s=\inf\{X_u,0\leq u\leq s\}.$$ 
By It\^o's excursion theory for
Markov processes, the excursions away from
$0$ of the process $X-\un{X}$ under $P$ are distributed according to a
Poisson point process that can be described by the It\^o excursion measure,
which we call
$N$. We now either consider the process
$X$ under the law
$P$ that makes it a L\'evy process starting at $0$, or under the 
$\sigma$-finite measure
$N$ under which the sample paths are excursions with finite lifetime
$\zeta$ (since $E[X_1]=0$).
Let $N^{(v)}$ be a regular version of the probability law $N(\cdot|\zeta=v)$, 
which is weakly continuous in $v$. That is, for any positive continuous
functional $G$, 
$$N(G)=\int_{(0,\infty)}N(\zeta\in \d v) N^{(v)}(G)$$ 
and $\lim
N^{(w)}(G)=N^{(v)}(G)$ as $w\to v$. Such a version can be obtained by
scaling: for any fixed $\eta>0$, the process 
$$\left((v/\zeta)^{1/\alpha}X_{\zeta s/v},0\leq
s\leq v\right)\quad\quad \mbox{ under }
N(\cdot|\zeta>\eta)=\frac{N(\cdot,\zeta>\eta)}{N(\zeta>\eta)}$$
%if $\un{g}_1$ (resp.\ $\un{d}_1$) denotes the  first time before
%(resp.\ after) $1$ when $X-\un{X}$ equals $0$,  then the law under $P$ of 
%$$X^{[\un{g}_1,\un{d}_1],t}_s=\frac{1}{(v(\un{d}_1-\un{g}_1))^{1/\alpha}}
%\left(X_{v(\un{d}_1-\un{g}_1)s}-\un{X}_{v(\un{d}_1-\un{g}_1)s}\right)\quad
%0\leq s\leq v$$
is $N^{(v)}$. See \cite{chaumont97} for this and other interesting ways to
obtain processes with law $N^{(v)}$ by path transformations. In  particular,
one has the scaling property at the level of  conditioned excursions: under
$N^{(v)}$, 
$\left(v^{-1/\alpha}X_{v s},0\leq s\leq 1\right)$
has law $N^{(1)}$. 

\paragraph{First-passage subordinator}
Let $T$ be the right-continuous inverse of the increasing process
$-\un{X}$, that is, 
$$T_x=\inf\{s\geq 0: \un{X}_s<-x\}.$$ 
Then it is 
known that under $P$, $T$ is a subordinator, that is, an increasing L\'evy
process. According to \cite[Theorem VII.1.1]{bertlev96}, its Laplace exponent
$\phi$ is the inverse function of the restriction of the Laplace exponent of
$X$ to $\R_+$. Thus $\phi(\lambda)=\lambda^{1/\alpha}$, and $T$ is a stable
subordinator with index $1/\alpha$, as defined above. The 
L\'evy-Khintchine formula gives,
$$E[\exp(-\lambda T_x)]=\exp(-x\lambda^{1/\alpha})=\exp\left(
x\int_0^{\infty}\frac{c_{\alpha}\d y}{y^{1+1/\alpha}}(1-e^{-\lambda
y})\right)\quad\mbox{ for } 
\quad \lambda,x\geq 0.$$
where $c_{\alpha}$ has been defined in the introduction.
Recall our assumption that $X$ has a marginal density at time $s$ 
under $P$, called $p_s(\cdot)$. Then under $P$, the inverse subordinator
$T$ has also bicontinuous densities, given e.g.\  by 
\cite[Corollary VII.1.3]{bertlev96}: 
\begin{equation}\label{ballotc}
q_x(s)=\frac{P(T_x\in \d s)}{\d s}=\frac{x}{s}p_s(x). 
\end{equation}
This equation can be derived from the ballot theorem of 
Tak\'acs \cite{takacs}. 

Let us now discuss the conditioned forms of distributions of the sequence
$\Delta T_{[0,x]}$. An easy way to obtain nice regular versions for  these
conditional laws is developed in
\cite{ppy92,pitmanpk02}, and uses the  notion of size-biased fragment.
Precisely, the range of any subordinator, with drift $0$ say (which we will
assume in the sequel), between times $0$ and
$x$, induces a partition of $[0,T_x]$ into subintervals with sum $T_x$. 
Consider a sequence $(U_i,i\geq 1)$ of independent uniform $(0,1)$ variables, 
independent of $T$, and let $\Delta^*_1(x),\Delta^*_2(x),\ldots$ be the 
sequence of the lengths of these intervals in the order in which they are 
discovered by the $U_i$'s. That is, $\Delta^*_1(x)$ is the length of the 
interval in which $T_xU_1$ falls, $\Delta^*_2(x)$ is the length of the 
first interval different from the one containing $T_xU_1$ 
in which $T_xU_i$ falls, and so on. Then Palm measure results for Poisson 
clouds give the following result (specialized to the case of stable
subordinators).

\begin{lmm}\label{permanform}
The joint law under $P$ of $(\Delta^*_1(x),T_x)$ is
\begin{equation}\label{sizeb}
P(\Delta^*_1(x)\in \d y, T_x\in \d s)=
\frac{c_{\alpha}xq_x(s-y)}{s y^{1/\alpha}}\d y \d s,
\end{equation}
and more generally for $j\geq 1$, 
$$P\left(\Delta^*_j(x)\in \d y\left|
T_x=s_0,\Delta^*_1(x)=s_1,\ldots,\Delta^*_{j-1}(x)=s_{j-1}
\right.\right)=\frac{c_{\alpha} x q_x(s-y)}{
sy^{1/\alpha}q_x(s)} \d y,$$
where $s=s_0-s_1-\ldots -s_{j-1}$.
\end{lmm}

This gives a nice regular conditional version for $(\Delta^*_i(x),i\geq 1)$
given $T_x$, and thus induces a conditional version for 
$\Delta T_{[0,x]}$ given $T_x$, by ranking, where 
$\Delta T_{[0,x]}$ is the sequence of jumps of $T$ before $x$, ranked in
decreasing order of magnitude. 

\subsection{The stable tree}\label{secsttree}

We now introduce the models of trees we will consider. 
This section is mainly inspired by \cite{duqleg02,duq02}. With the notations 
of section \ref{levyproc}, for 
$u\geq 0$, let $R^{(u)}$ be the time-reversed process of $X$ at time $u$:
$$R^{(u)}_s=X_u-X_{(u-s)-}\quad ,\quad 0\leq s\leq u.$$
It is standard that this process has the same law as $X$ killed
at time $u$ under $P$. Let also 
$$\ov{R}^{(u)}_s=\sup_{0\leq v\leq s} R^{(u)}_v \quad,  \quad 0\leq s\leq u$$
be its supremum process. 
We let $H_u$ be the local time at $0$ 
of the process $R^{(u)}$ reflected under its supremum 
$\ov{R}^{(u)}$ up to time $u$. The normalization can be chosen so that 
$$H_u=\lim_{\eps\downarrow0}\frac{1}{\eps}\int_0^u\ind_{\{
\ov{R}^{(u)}_s-R^{(u)}_s\leq \eps\}}\d s$$
It is known by \cite[Theorem 1.4.3]{duqleg02} that $H$ admits a continuous
version, with which we shall work in the sequel. 
It has to be noticed that $H$ is not a Markov process (the only exception
in the theory of L\'evy trees is the Brownian tree obtained when $P$ is
the law of Brownian motion with drift, which has been excluded in our
discussion). As a matter of fact, it can be checked that
$H$ admits local minima that are attained an infinite number of times as
soon as
$X$ has jumps, a property that sounds strange by contrast with most of the
usually studied  stochastic processes. To see this, consider a jump time $t$
of $X$, and  let $t_1,t_2>t$ so that $\inf_{t\leq u\leq t_i}X_u=X_{t_i}$ and 
$X_{t-}<X_{t_i}<X_t$, $i\in\{1,2\}$. Then it is easy to see that 
$H_{t}=H_{t_1}=H_{t_2}$ and that one may in fact find an infinite number
of distinct $t_i$'s satisfying the properties of $t_1,t_2$. On the other hand, 
it is not difficult to see that $H_t$ is a local minimum of $H$. 

It is  shown in \cite{duqleg02} that the definition of $H$ 
still makes sense under the 
$\sigma$-finite measure $N$ rather than the probability law $P$. The process
$H$ is then defined only on $[0,\zeta]$, and we call it the {\em excursion
of the height process}.
One can define without difficulty, using the scaling property, 
the height process under the laws $N^{(v)}$: this is simply the law of 
$$\left(\left(\frac{v}{\zeta}\right)^{1-1/\alpha}
H_{\zeta  t/v},0\leq t\leq v\right)
\quad \mbox{ under } N(\cdot,\zeta>\eta)$$
Call it the law of the 
excursion of the height process with duration $v$. The following scaling
property is the key for the self-similarity of $F^-$: for every $x>0$, 
\begin{equation}\label{scalep}
(v^{1/\alpha-1}H_{sv},0\leq s\leq 1) \mbox{ under }N^{(v)} 
\build=_{}^{d}(H_s,0\leq s\leq 1) \mbox{ under }N^{(1)}.
\end{equation}
This property is inherited from the scaling property of $X$, and it is 
easily obtained e.g.\ by the above definition of $H$ as an approximation.

%\noindent{\bf Remark. }
%We believe that it may 
%be also possible to construct such conditioned excursions of the height 
%process by adapting the arguments of Duquesne \cite{duq02} in the general 
%case, by first defining the local time of the dual processes for bridges
%(this could be done using the absolute continuity
%relation used to define Markovian bridges) 
%and then using Vervaat's theorem. We could thus
%define the excursion of the 
%height process with duration $v$ for any $v>0$, as soon as we
%could prove a certain continuity relation between these different
%conditioned excursions of the height process. But we will avoid such
%discussions in the  sequel. 

An important tool for studying the height process is its
{\em local time process}, or {\em width process}, which we will denote by 
$(L^t_s,t\geq 0,s\geq0)$. It can be obtained a.s.\ for every fixed $s,t$ by 
$$L^t_s=\lim_{\eps\downarrow0}
\frac{1}{\eps}\int_0^s\ind_{\{t< H_u\leq t+\eps\}}\d u.$$
$L_s^t$ is then the density of the occupation measure of $H$ at level
$t$ and time $s$. For $t=0$,
one has that $(L_s^0,s\geq 0)$ is the inverse of the subordinator 
$T$, which is a reminiscent of the fact that the excursions of the height 
process are in one-to-one correspondence with excursions of $X$ with the 
same lengths.
According to the Ray-Knight theorem \cite[Theorem 1.4.1]{duqleg02}, for every
$x>0$, the process $(L^t_{T_x},t\geq 0)$ is a continuous-time branching 
process with branching mechanism $\lambda^{\alpha}$, 
in short $\alpha$-CSBP. We will
recall basic and less basic features about this processes in Sect.\ 
\ref{secCSBP}, where in particular an interpretation for the law of the 
process $(L^t_1,t\geq 0)$ under $N^{(1)}$ will be given. 
%The key feature of this process is:
%\begin{prp}\label{rkthm}{\rm \cite[Theorem
%1.4.1]{duqleg02}}
%For any $x>0$, the process $(L_{T_x}^t,t\geq0)$ is the $\psi$-continuous
%state branching process started at $x$ (in short the $\psi$-CSBP), its 
%law $\P_x$ is the unique probability measure that makes the process 
%$(X_t,t\geq0)$ a right-continuous Markov process
%starting at $x$ with transition probabilities characterized by 
%$$\E[\exp(-\lambda X_{t+r})|X_t=y]=\exp(-yu_r(\lambda)),$$
%where $u_r(\lambda)$ is determined by the equation
%$$\int_{u_r(\lambda)}^{\lambda}\frac{\d u}{\psi(u)}=r.$$
%\end{prp}
For now we just note that for every $x$ the process $(L_{T_x}^t,t\geq0)$ is a 
process with no negative jumps, and a jump of this process at time
$t$ corresponds 
precisely to one of the infinitely often attained local infima of the 
height process. With the forthcoming interpretation of the tree encoded within
excursions of the height process, this means that there is a branchpoint
with infinite degree at level $t$. It is again possible to define
the local time process under the excursion measure $N$, and by scaling
it is also possible
to define the local time process under $N^{(v)}$. 
%From the 
%last proposition, one can get interpretations for the process $(L_{\zeta}^t,
%t\geq 0)$ under the measure $N$ or of $(L_v^t,t\geq 0)$ under $N^{(v)}$ in
%terms of conditioned continuous-state branching processes (branching processes
%conditioned to have a fixed total population and conditioned to leave 
%from the exit value $0$). These in turn have been studied by Lambert 
%\cite{lambert01} in terms
%of immigration. 

Let us now motivate the term of ``height process'' for $H$. 
Under the $\sigma$-finite ``law'' $N$, we define a tree structure following 
\cite{aldouscrt93,legall93}. 

First we introduce some extra vocabulary. Let ${\bf T}$ be the set of 
finite rooted plane trees, that is, for any 
$\TT\in{\bf T}$, each set of children of 
a vertex $v\in \TT$ is ordered as first, second, ..., last child. Let 
${\bf T}^*\subset {\bf T}$ be those rooted plane trees for which the 
out-degree (number of children) of vertices is never $1$. Let 
${\bf T}_n$ and ${\bf T}^*_n$ be the corresponding sets of trees that have 
exactly $n$ leaves (vertices with out-degree $0$).
A marked tree $\vartheta$ is a pair $(\TT,\{h_v,v\in\TT\})$ where 
$\TT\in{\bf T}$ and $h_v\geq 0$ for every vertex $v$ of $\TT$ (which we 
denote by $v\in\TT$). The tree $\TT$ is called the skeleton of $\vartheta$, 
and the $h_v$'s are the marks. These marks induce a distance  
tree, given by $d_{\vartheta}(v,v')=\sum_{w\in[[v,v']]} h_w$
if $v,v'\in \vartheta$ are two vertices
of the marked tree,
where $[[v,v']]$ is the set of vertices of the 
path from $v$ to $v'$ in the skeleton. 
The distance of a vertex to the root will be called its {\em height}.
Let $\mathbb{T}^*_n$ be the set of 
marked trees with $n$ leaves and no out-degree equal to $1$. 

Let $(U_i,i\geq 1)$ be independent random variables with uniform law on 
$(0,1)$ and independent of the excursion $H$ of the height process. 
One may define a random marked tree $\vartheta(U_1,\ldots,U_k)=\vartheta_k\in
\mathbb{T}^*_k$, as follows. For $u,v\in[0,\zeta]$ let 
$m(u,v)=\inf_{s\in[u,v]}H_s$. Roughly, the key fact about $\vartheta_k$
is that the height 
of the $i$-th leaf to the 
root is $H_{U_{(i)}}$, 
where $(U_{(i)},1\leq i\leq k)$ are the order statistics 
of $(U_i,1\leq i\leq k)$, and the ancestor of the $i$-th and $j$-th leaves
has height $m(\zeta U_{(i)},\zeta U_{(j)})$ for every $i,j$. This allows 
to build recursively a tree by first putting the mark
$h_{\rm root}=\inf_{1\leq i<j\leq k} m(U_i,U_j)$ on a root vertex.
Let $c_{\rm root}$ 
be the number of excursions of $H$ above level $h_{\rm root}$ 
in which at least one 
$\zeta U_i$ falls. Attach $c_{\rm root}$ vertices to the root, and let the 
$i$-th of these vertices be the root of the tree embedded in the $i$-th
of these excursions above level $h_{\rm root}$. Go on until the excursions
separate the variables $U_i$. By construction $\vartheta_k\in\mathbb{T}^*_k$.
Adding a $(k+1)$-th variable $U_{k+1}$ to the first $k$ just adds a new branch 
to the tree in a consistent way as $k$ varies. 
%We call the family $(\vartheta_k,k\geq 1)$ of marked trees
%the  tree. One may think of this globally by letting $k\to\infty$ and
%completing the metric space induced by the variables. 
%In this case one may define
%a mass measure inherited from Lebesgue measure 
%on $[0,\zeta]$ and supported by the leaves, and a length measure induced 
%by the distance
%$$d(u,v)=H_u+H_v-2m(u,v),\quad u,v\in[0,\zeta].$$
%More generally, the height process $(H_t,t\geq 0)$ can be interpreted as
%a forest of continuum random trees, each excursion of $H$ away from
%$0$ corresponding to a tree of the forest. 

As noted above, we may as well define the trees $(\vartheta_k,k\geq 0)$ under
the law $N^{(1)}$ by means of scaling. 
\begin{defn}\label{sttre}
The family of marked trees $(\vartheta_k,k\geq 1)$ associated with 
the height process under the law $N^{(1)}$ is called the stable tree.
\end{defn}

\noindent{\bf Remark. }
The previous definition is not the only way to characterize the 
same object. After all, we could have called the height process $H$ under 
$N^{(1)}$ itself the stable tree. 
Alternatively, one easily sees that the marked tree 
$\vartheta_k$ can be interpreted as a subset of $l^1$, each new branch
going in a direction orthogonal to the preceding branches, in a consistent 
way as $k$ varies. Then it makes sense to take the completion of
$\cup_{k\geq 1}\vartheta_k$, which we could also call the stable 
tree. The distance on the tree then corresponds to the 
metric defined under $N^{(1)}$ by 
$$d(u,v)=H_u+H_v-2m(u,v),\quad u,v\in[0,1].$$
With this way of looking at things, the leaves of the tree are 
uncountable and everywhere dense in the tree, and the empirical distribution 
on the leaves of $\vartheta_k$ converges weakly to a probability 
measure on 
the stable tree, called the mass measure. Then it turns out that 
$\vartheta_k$ 
is equal in law to the subtree of the stable tree that is 
spanned by the root and $k$ independent leaves distributed according to 
the mass measure. Hence, the mass measure 
is represented by Lebesgue measure on $[0,1]$ in the 
coding of the stable tree through its height process. This is coherent with 
the definition of $F^-(t)$ as the ``masses of the tree components located
above height $t$''. The equivalence between these possible definitions 
is discussed in \cite{aldouscrt93}.

The key property for obtaining the dislocation measure of $F^-$
is the following description of the law of the skeleton of
$\vartheta_n$, and  the mark of the root of $\vartheta_1$.
For $\TT\in{\bf T}$ let 
${\cal N}_{\TT}$ be the set of non-leaf vertices of $\TT$ and for $v\in
{\cal N}_{\TT}$ let
$c_v(\TT)$ be the number of children of $v$.
From the more complete description of the marked trees in 
\cite[Theorem 3.3.3]{duqleg02}, we recall that

\begin{prp}\label{skelstruc}
The probability that the skeleton of $\vartheta_n$ is 
$\TT\in {\bf T}^*_k$ is 
$$\frac{n!}{(\alpha-1)(2\alpha-1)\ldots((n-1)\alpha-1)}\prod_{v\in 
{\cal N}_{\TT}}\frac{|(\alpha-1)(\alpha-2)\ldots(\alpha-c_v(\TT)+1)|}{
c_v(\TT)!}.$$
Moreover, the law of the mark of the root in $\vartheta_1$ is 
$$N^{(1)}(H_{U_1}\in \d h)=\alpha \Gamma\left(1-\frac{1}{\alpha}\right)
\chi_{\alpha h}(1)\d h,$$ 
where $\chi_x(s)$ is the 
density of the stable $1-1/\alpha$ subordinator (with Laplace exponent
equal to $\lambda^{1-1/\alpha}$) at time $x$.
\end{prp}

\subsection{Some results on self-similar fragmentations}\label{ssfrag}

In this section we are going to recall some basic facts about the theory
of self-similar fragmentations, and also introduce some useful 
ways to recover the characteristics of these fragmentations. We will 
suppose that the fragmentations we consider are not trivial, that is, 
they are not equal to their initial state for every time.
It will be useful to consider not only $S$-valued (or {\em ranked}) 
fragmentations, but also fragmentations with values in the set 
of open subsets of $(0,1)$ and in the set of partitions of
$\N=\{1,2,\ldots\}$, respectively called {\em interval} and {\em
partition-valued} fragmentations. As established in
\cite{bertsfrag02,berest02},  there is a one-to-one mapping between the laws
of the  three kinds of fragmentation when they satisfy a self-similarity
property  that is similar to that of the ranked fragmentations. That is, 
each of them is characterized by the same 3-tuple $(\beta,c,\nu)$ introduced 
above. To be completely accurate, we should stress that there actually exist
several versions of interval partitions that give the same ranked or 
partition-valued fragmentation, but all these versions have the 
same characteristics $(\beta,c,\nu)$. Let us make the terms precise.

Let ${\cal P}$ be the set of unordered partitions of $\N$. 
An exchangeable partition $\Pi$ is a ${\cal P}$-valued random variable
whose restriction $\Pi_n$ to $[n]=\{1,\ldots,n\}$ has an invariant law 
under the action of the permutations of $[n]$. By Kingman's representation
theorem \cite{kingman78,aldous85}, the blocks of 
exchangeable partitions of $\N$ admit almost-sure asymptotic frequencies, 
that is, if $\Pi=\{B_1,B_2,\ldots\}$ where the $B_i$'s are listed by order
of their least element, then 
$$\Lambda(B_i)=\lim_{n\to\infty}\frac{|B_i\cap[n]|}{n}$$
exists a.s.\ for every $i\geq 0$. Denoting by $\Lambda(\Pi)$ the ranked 
sequence of these asymptotic frequencies, $\Lambda(\Pi)$ is then a 
$S$-valued random variable, whose law characterizes that of $\Pi$.

A self-similar
partition-valued fragmentation $(\Pi(t),t\geq 0)$ with index $\beta$
is a ${\cal P}$-valued 
c\`adl\`ag process that is continuous in probability, exchangeable, 
meaning that for every permutation $\sigma$ of $\N$, 
$(\sigma \Pi(t),t\geq 0)$ and 
$(\Pi(t),t\geq 0)$ have the same law, and such that given
$\Pi(t)=\{B_1,B_2,\ldots\}$, the variable $\Pi(t+t')$ has the law of 
the partition with blocks $\Pi^{(i)}(\Lambda(B_i)^{\beta}t')
\circ B_i$ where the $\Pi^{(i)}$ are independent copies of $\Pi$. Here , the 
operation $\circ$ is the natural ``fragmentation'' operation of a set by a 
partition: if $\Pi=\{B_1,B_2,\ldots\}$ and $C\subset\N$, then 
$\Pi\circ C$ is the partition of $C$ with blocks $B_i\cap C$.

A self-similar interval partition $(I(t),t\geq 0)$ with index $\beta$
is a process with 
values in the open subsets of $(0,1)$ which is right-continuous and
continuous in probability for the usual Hausdorff metric, with the property 
that given $I(t)=\cup_{i\geq 1}I_i$ say, where the $I_i$ are the disjoint 
connected components of $I(t)$, the set $I(t+t')$ has the law of 
$\cup_{i\geq 1}g_i(I^{(i)}(t' |I_i|^{\beta}))$, where $|I_i|$ is the 
length of $I_i$, $g_i$ is the affine transformation that maps $(0,1)$ to 
$I_i$ and conserves orientation and the $I^{(i)}$ are independent copies of 
$I$. 

Consider an interval self-similar fragmentation $(I(t),t\geq 0)$, with
characteristic $3$-tuple $(\beta,0,\nu)$ (the case when $c>0$ would be similar, but we
do not need it in the sequel).  Let $U_i,i\geq 1$ be independent uniform random
variables on $(0,1)$. These induce a partition-valued fragmentation
$(\Pi(t),t\geq 0)$ by saying that $i\build\sim_{}^{\Pi(a)}j$ iff $U_i$ and
$U_j$ are in the same  connected component of $I(t)$. It is known
\cite{bertsfrag02} that $\Pi$ is a self-similar fragmentation with values in
the set of partitions of $\N$ and characteristics $(\beta,0,\nu)$. For $n\geq
2$ let ${\cal P}_n^*$ be the set of partitions of $\N$ whose restriction to
$[n]$ is non-trivial, i.e.\ different from $\{[n]\}$. Then there is some random
time $t_n>0$ such that the restriction of $\Pi(t)$ to $[n]$ jumps from the
trivial state $\{[n]\}$ to some non-trivial  state at time $t_n$. Let $\rho(n)$
be the law of the restriction of $\Pi(t_n)$ to $[n]$. The next Lemma states
that the knowledge of the family $(\rho(n),n\geq 2)$ almost determines the
dislocation measure $\nu$ of the fragmentation. Precisely, we introduce from
\cite{berthfrag01} the notion of characteristic measure of the fragmentation.
This measure, denoted by $\kappa$, is a $\sigma$-finite measure supported by
the non-trivial partitions of $\N$, which is determined by the dislocation
measure of the fragmentation. Precisely, this measure may be written as 
$$\kappa(\d\pi)=\int_S\nu(\d\bs)\kappa_{\bs}(\d\pi),$$
where $\kappa_{\bs}$ is the law of the exchangeable partition of $\N$ with
ranked asymptotic frequencies given by $\bs$. Conversely, this measure
characterizes the dislocation measure $\nu$ (simply by taking the asymptotic
frequencies of the typical partition under $\kappa$). 

\begin{lmm}\label{rhocar}
The restriction of $\kappa$ to the non-trivial partitions of $[n]$, for $n\geq
2$, equals $q(n)\rho(n)$, for some sequence $(q(n),n\geq 2)$ of strictly
positive numbers. As a consequence, the dislocation measure of the
fragmentation $I$ is characterized by the sequence of laws $(\rho(n),n\geq 2)$,
up to a multiplicative constant. 
\end{lmm}

Otherwise said, and using the correspondence between self-similar
fragmentations with same dislocation measure and different indices
established by Bertoin \cite{bertsfrag02} by introducing the
appropriate time-changes, if we have two interval-valued self-similar
fragmentations $I$ and $I'$ with the same index and no erosion, and with the
same associated probabilities $\rho(n)$ and $\rho'(n)$, $n\geq 1$, then there
exists $K>0$ such that $(I(Kt),t\geq 0)$ has the same dislocation measure as
$I'$.

\begin{proof}
Suppose $\beta=0$, then the result is almost immediate by the results
of \cite{berthfrag01} on homogeneous fragmentation processes. In this case
$q(n)$ is the inverse of the expected jump time
of $\Pi$ in ${\cal P}_n^*$, and the restriction of the measure
$q(n+1)\rho(n+1)$ to the set of non-trivial partitions of $[n]$ is
$q(n)\rho(n)$, for every $n\geq 1$. Hence, it is easy to see that the knowledge
on $\rho(n)$ determines uniquely the sequence $(q(n),n\geq 1)$, up to  a
multiplicative positive constant: one simply has
$q(n)/q(n+1)=\rho(n+1)(\pi|_{[n]}:\pi\in{\cal P}_n^*)$, where $\pi|_{[n]}$
denotes the restriction of $\pi$ to $[n]$. It remains to notice that the
sequence of restrictions $(q(n)\rho(n),n\geq 2)$ characterizes $\kappa$. 
%Now we know that there
%exists a
%$\sigma$-finite measure
%$\kappa$ on the set of partitions of $\N$ that puts no mass on the trivial
%partition, and such that the restriction of $\kappa$ to the non-trivial
%partitions of $[n]$ is
%$q(n)\rho(n)$. This measure is known to characterize the law of the
%fragmentation, and in particular $\nu$, which is obtained as the law of the
%almost-surely existing ranked limiting frequencies of the generic partition
%under $\kappa$. 
%Hence, given the family
%$(\rho(n),n\geq 1)$ one deduces the exchangeable dislocation measure
%$\kappa$ of $\Pi$ up to a positive multiplicative constant.

When $\beta\neq0$, we obtain the same results by noticing that the law
$\rho(n)$ still equals the law of the restriction to $[n]$ of the exchangeable
partition with limiting frequencies having the ``law'' $\nu$ and restricted to
${\cal P}_n^*$, up to a multiplicative constant. Indeed, let $I^*(t)$ be the
subinterval of $I(t)$ containing $U_1$ at time $t$, and recall
\cite{bertsfrag02} that if 
$$a(t)=\inf\left\{u\geq 0:\int_0^u|I^*(v)|^{\beta}\d v>t\right\},$$
then $(|I^*(a(t))|,t\geq 0)$ evolves as the fragment containing $U_1$ in an
interval fragmentation with characteristics $(0,0,\nu)$. Now, before time
$t_n$, the fragment containing $U_1$ is the same as that containing all the
$(U_i,1\leq i\leq n)$. Hence, $a(t_n)$ is the first time when $\Pi'$ jumps
in ${\cal P}_n^*$ for some homogeneous partition-valued fragmentation process
$\Pi'$ with characteristics $(0,0,\nu)$, and the law of $\Pi'(a(t_n))$
restricted to $[n]$ is $\rho(n)$. Hence the result. \cq
\end{proof}

We also cite the following result \cite[Proposition 3]{MS03} which allows to 
recover the dislocation measure of a self-similar fragmentation with 
positive index out of its semigroup. We will not use this proposition 
in a proof, but it is useful to keep it in mind to conjecture the 
form of the dislocation measure of $F^-$, as it will be done below. 

\begin{prp}\label{msdisl}
Let $(F(t),t\geq 0)$ be a ranked self-similar fragmentation with
characteristics $(\beta,0,\nu)$, $\beta\geq 0$. Then for every continuous
bounded function $G$ on $S$ which is null on an open neighborhood of
$(1,0,\ldots)$, one has
$$\frac{1}{t}E[G(F(t))]\build\to_{t\downarrow0}^{}\nu(G).$$
\end{prp}

\section{Study of $F^-$}\label{fmoins}

We now specifically turn to the study of $F^-$ defined in the introduction.
Although some of the  results below may be easily generalized to a broader
``L\'evy context'', we will suppose in this section that $X$ is a stable
process with  index $\alpha\in(1,2)$, with first-passage subordinator $T$. The
index 
$\alpha$ will be dropped from the notation by contrast with 
the introduction. The references to 
height processes, excursion measures and so on, will always be
with respect to this process, unless otherwise specified. Also, for the needs
of the proofs below, we define the process $(F^-(t),t\geq 0)$ not only
under the law $N^{(1)}$ used to define the stable tree, but also for all the
excursion measures $N^{(v)}$ and $N$. Under $N^{(v)}$, let $F^-(t)$ be the
decreasing sequence of lengths of the constancy intervals of
$I_-(t)=\{s\in(0,v):H_s>t\}$ ($v$ is replaced by $\zeta$ under $N$). To avoid
confusions, we will always mention in Sect.\ \ref{selfsemig} the measure we are
working with, but this formalism will be abandoned in the following sections
where no more use of
$N^{(v)}$ is made with $v\neq1$. 

The study contains four steps.
First we prove that self-similarity property for
$F^-$ and  make its semigroup explicit. Heuristic
arguments based on generators of conditioned CSBP's allow us to
conjecture the rough shape of the dislocation measure . Then we  prove that
the erosion coefficient is
$0$ by studying the evolution of a tagged fragment. We are then able to apply
Lemma \ref{rhocar}, giving us the dislocation measure up to a constant, and we
finally recover the constant by re-obtaining the results needed in the
second step by another computation. 

\subsection{Self-similarity and semigroup}\label{selfsemig}

The self-similarity and the description of the semigroup rely strongly on the 
following result, which is a variant of 
\cite[Proposition 1.3.1]{duqleg02}. 
For $t,s\geq 0$ let 
$$\gamma^t_s=\inf\{u\geq 0: \int_0^u \ind_{\{H_v>t\}}\d v >s\}$$ 
and 
$$\widetilde{\gamma}^t_s=\inf\{u\geq 0:\int_0^u\ind_{\{H_v\leq
t\}}\d v>s\}.$$   
Denote by ${\cal H}_t$ the sigma-field generated by the
process 
$(H_{\widetilde{\gamma}^t_s},s\geq 0)$ and the $P$-negligible sets. Let
also
$(H^t_s,s\geq 0)$ be the process $(H_{\gamma^t_s}-t,s\geq 0)$. Then under 
$P$, $H^t$ is 
independent of ${\cal H}_t$, and its law is the same as that of 
$H$ under $P$.

As a first consequence, we obtain that the excursions of $H$ above level 
$t$, that is, the excursions of $H^t$ above level $0$, are, conditionally on 
their durations, independent excursions of $H$. This simple result allows
us to state the Markov property and self-similarity of $F^-$. 
%Notice that if we denote by $L^t_{T_1}$ the local time
%of $H$ at level $t$ and up to time $T_1$, then this is plainly equal to 
%the local time of $H^t$ at time $\gamma^t_{T_1}$ and level $0$. We claim that
%this is in fact independent of $H^t$. To see this, notice that this local
%time may be also obtained as the occupation density of 
%$(H_{\widetilde{\gamma}^t_s},s\geq 0)$ at level $t$ and time 
%$\widetilde{\gamma}^t_{T_1}$. As such, 
%this is a ${\cal H}_t$-measurable variable, hence independent of $H^t$. Now, 
%conditionally on $L^t_{T_1}=\ell$, 
%we have that $F^-_1(t)$ is the sequence of 
%lengths of the constancy intervals of $\{s\in[0,T^t_{\ell}]:H^t_s>0\}$ where 
%$T^{t}$ is the inverse local time of $H^t$ at level $0$, and this has the 
%same law as $T$. Thus, this is exactly the law of lengths of excursions 
%of $X$ above its infimum process, up to time $T_{\ell}$. As a consequence, 
%\begin{equation}
%P(F^-_1(t)\in \d\bs)=\int_0^{\infty}P(L^t_{T_1}\in \d \ell)P\left(
%\Delta T_{[0,\ell]}\in \d\bs\right).
%\label{sf1}
%\end{equation}
%From this we deduce the self-similarity property and semigroup 
In the following statement, it has to be understood that 
we work under the probability $N^{(1)}$ and that the process $H$ that is 
considered is the same that is used to construct $F^-$. 

\begin{lmm}\label{spropfm}
Conditionally on $F^-(t)=(x_1,x_2,\ldots)$, the excursions of $H$ above level
$t$, that is, of $H^t$ away from $0$, are independent excursions with
respective laws $N^{(x_1)},N^{(x_2)},\ldots$. 

As a consequence, the process $F^-$ is a self-similar fragmentation process
with index $1/\alpha -1$. 
\end{lmm}

\begin{proof}
By the previous considerations on $H^t$, we have that under $P$, given 
that the lengths of interval components of the set $\{s\in[0,T_1]:H_s>t\}$
ranked in decreasing order are equal to
$(x_1,x_2,\ldots)$, the excursions of the killed process 
$(H(t),0\leq t\leq T_1)$ above level $t$ are independent excursions of 
$H$ with durations $x_1,x_2,\ldots$. The first part of the statement follows
by  considering the first excursion of $H$ (or of $X$)
that has duration greater than some $v>0$, which gives the result under the 
measure $N(\cdot,\zeta>v)$, hence for $N$, hence for $N^{(v)}$ for almost 
all $v$, and then for $v=1$ by continuity of the measures $N^{(v)}$. 
%Hence conditionally on $F^-(t)=(x_1,x_2,\ldots)$, the process $F(t+t'),t'\geq
%0)$ follows the law of the decreasing rearrangement of the processes 
%$(x_i F^-_i(t'),t'\geq 0)$

Thus, conditionally on $F^-(t)=(x_1,x_2,\ldots)$, 
the process $(F^-(t+t'),t\geq 0)$ has the same
law as the random sequence 
obtained by taking independent excursions $H^{(x_1)},H^{(x_2)},
\ldots$ with durations
$x_1,x_2,\ldots$ of the height process, and then arranging in decreasing order 
the lengths of constancy intervals of the sets
$$\{s\in[0,x_i]:H^{(x_i)}_s>t'\}.$$
It thus follows from the scaling property (\ref{scalep}) 
of the excursions of $H$ 
that conditionally on $F^-(t)=(x_1,x_2,\ldots)$, the
process $(F^-(t+t'),t'\geq 0)$ has the same law as the decreasing
rearrangement of the processes $(x_i F^-_{(i)}(x_i^{1/\alpha-1}t'),t'\geq 0)$,
where the $F_{(i)}^-$'s are independent copies of $F^-$. The fact that
$F^-$ is a Markov process that is continuous in probability easily follows, 
as does the self-similar fragmentation property with the index
$1/\alpha-1$. 
\cq
\end{proof}

We now turn our attention to the semigroup of $F^-$.

\begin{prp}\label{semigrmoins}
For every $t\geq 0$ one has
\begin{equation}\label{semifm}
N^{(1)}(F^-(t)\in \d\bs)=\int_{\R_+\times [0,1]}N^{(1)}\left(L_1^t\in \d \ell,
\int_t^{\infty}
\!\!\d b\, L_1^b\in \d z\right)P\left(\Delta T_{[0,\ell]}\in \d\bs
\left|T_{\ell}=z\right.\right),
\end{equation}
with the convention that 
the law $P(\Delta T_{[0,0]}\in \d\bs|T_0=z)$ is the Dirac mass at the sequence 
$(z,0,0\ldots)$ for every $z\geq 0$. 
\end{prp}

\begin{proof} 
It suffices to prove the result for some fixed $t>0$. 
Let $\omega(t)=\inf\{s\geq 0:H_s>t\}$,
$d_{\omega(t)}=\inf\{s\geq \omega(t):X_s=\un{X}_s\}$
and $g_{\omega(t)}=\sup\{s\leq \omega(t):X_s=\un{X}_s\}$.  
Call $\FF^-(t)$ the ranked sequence of the lengths of the interval 
components of the set $\{s\in[\omega(t),d_{\omega(t)}]:H_s>t\}$. Notice that 
under the law $N^{(1)}$, $\FF^-$ would be $F^-$, but we will first define
$\FF^-$ under $P$.  By the definition of $H$, $\omega(t)$ and $d_{\omega(t)}$
are stopping times with respect to the natural filtration generated by $X$. In
fact, it also holds that $\omega(t)$ is a {\em terminal time}, that is, 
$$\omega(t)=s+\inf\{u\geq 0: H_{s+u}>t\} \quad\mbox{ on } \{\omega(t)>s\}.$$
Moreover, $0<\omega(t)<\infty$ $P$-a.s., 
because of the continuity of $H$ and the fact that excursions of $H$ have a
positive probability to hit level $t$ (which follows e.g. by scaling).  Recall
the notations at the beginning of the section,  and denote by $A^t$ and
$\wt{A}^t$ the  right-continuous inverses of $\gamma^t$ and $\wt{\gamma}^t$.
Then  the local time $L_{d_{\omega(t)}}^t$ is the local time at level $0$ and
time 
$A^t_{d_{\omega(t)}}$ of the process $H^t$. This is also equal to 
the local time of $(H_{\wt{\gamma}^t_s},s\geq 0)$ at level $t$ and time 
$\wt{A}^t_{d_{\omega(t)}}$. 
This last time is ${\cal H}_t$-measurable, as it is the first
time the process $(H_{\wt{\gamma}^t_s},s\geq 0)$ hits back $0$ after first 
hitting $t$. Hence $L^t_{d_{\omega(t)}}$ is ${\cal H}_t$-measurable, 
hence independent of $H^t$. Let $T^t$ be the inverse local time of 
$H^t$ at level $0$, which is $\sigma(H^t)$-measurable, hence independent of 
${\cal H}_t$, and has same law as $T$ since $H^t$ has same law as 
$H$ under $P$. Notice that 
$\FF^-(t)$ equals the sequence
$\Delta T^t_{[0,L^t_{d_{\omega(t)}}]}$, and that the $\sigma(H^t)$-measurable
random variable 
$\int_t^{\infty}\d b\, L^b_{d_{\omega(t)}}=T^t(L^t_{d_{\omega(t)}})$.
Thus, conditionally on $L^t_{d_{\omega(t)}}=\ell$ and 
$\int_t^{\infty}\d b\, L^b_{d_{\omega(t)}}=z$, 
$\FF^-(t)$ has law $P(\Delta T_{[0,\ell]}\in \d\bs|T_{\ell}=z)$. Hence  
$$P(\FF^-(t) \in \d\bs)=\int_{\R_+\times\R_+}P\left(L^t_{d_{\omega(t)}}\in 
\d\ell\, ,\, \int_t^{\infty}\!\! \d b\, L^b_{d_{\omega(t)}}\in \d z\right)
P(\Delta T_{[0,\ell]}\in \d\bs|T_{\ell}=z),$$
and also, since $d_{\omega(t)}-g_{\omega(t)}=\int_0^{\infty}\d b\, 
(L^b_{d_{\omega(t)}}-L^b_{g_{\omega(t)}})$ and since $\int_0^t\d b\,
(L^b_{d_{\omega(t)}}-L^b_{g_{\omega(t)}})$ is  independent of $\sigma(H^t)$,
the result also holds conditionally on 
$d_{\omega(t)}-g_{\omega(t)}$, namely
\begin{eqnarray*}
\lefteqn{P(\FF^-(t) \in \d\bs|d_{\omega(t)}-g_{\omega(t)})}\\
&=&
\int_{\R_+\times\R_+}P\left(L^t_{d_{\omega(t)}}\in 
\d\ell\, ,\, \int_t^{\infty}\!\! \d b\, L^b_{d_{\omega(t)}}\in \d z\bigg|
d_{\omega(t)}-g_{\omega(t)}\right)
P(\Delta T_{[0,\ell]}\in \d\bs|T_{\ell}=z).
\end{eqnarray*}
Now notice that the excursion of $H$ straddling time $\omega(t)$ is the first 
excursion of $H$ that attains level $t$, and apply 
\cite[Proposition XII.3.5]{revyor} to obtain that 
$$P(\FF^-(t) \in \d\bs|d_{\omega(t)}-g_{\omega(t)}=v)
=N^{(v)}(\zeta>\omega(t))^{-1}\, N^{(v)}(F^-_1(t) \in \d\bs,\, v>\omega(t)),$$
and similarly
\begin{eqnarray*}
\lefteqn{
P\left(L^t_{d_{\omega(t)}}\in 
\d\ell\, ,\, \int_t^{\infty}\!\! \d b\, L^b_{d_{\omega(t)}}\in \d z\bigg|
d_{\omega(t)}-g_{\omega(t)}=v\right)}\\
&=&
N^{(v)}(\zeta>\omega(t))^{-1}\, N^{(v)}\left(L^t_v\in 
\d\ell\, ,\, \int_t^{\infty}\!\! \d b\, L^b_v\in \d z\, ,\, v>
\omega(t)\right),
\end{eqnarray*}
for almost every $v$. This is generalized for every $v$ by continuity of the 
family $N^{(v)}$. 
Finally, notice that $\FF^-(t)=F^-(t)$ under $N$ and the $N^{(v)}$'s and that
we may remove the  indicator of $v>\omega(t)$ since
a.s.\ under $N^{(v)}$, $L^t_v=0$ if and only if $\max H\leq t$, to obtain 
$$N^{(v)}(F^-(t)\in \d\bs)=N^{(v)}\left(L_v^t\in \d \ell,
\int_t^{\infty}\d b\,
L_v^b\in \d z\right)P(\Delta T_{[0,\ell]}\in \d\bs|T_{\ell}=z).$$
Taking $v=1$ entails the claim.
\cq
\end{proof}

%As a consequence of this we obtain that the dislocation measure of $F^-$ does
%not charge the set $\{s\in S:\sum_{i\geq 1}s_i<1\}$. Indeed, the total mass of
%the fragments at time $t$ has the law of $\int_t^{\infty}\d s\, L_1^s$ under
%$N^{(1)}$. Since $(L_1^t,t\geq 0)$ is right continuous, the integral above is
%continuous as $t$ varies, so there is no sudden loss of mass. We will
%re-obtain this in the next section anyway. But more 
As a consequence of this result we may
conjecture the shape of the dislocation measure of $F^-$. The next subsections
will give essentially the rigorous proof of this conjecture, but
finding $\nu_-$ directly from the forthcoming computations would certainly
have been tricky without any former intuition. Roughly, suppose that the
statement of Proposition 
\ref{msdisl} remains true for negative self-similarity indices (which is
probably true, but we will not need it anyway). Then take $G$ a bounded
continuous function that is null on a neighborhood of $(1,0,\ldots)$ and write
$$N^{(1)}(G(F^-(t)))=\int_{\R_+\times [0,1]}N^{(1)}\left(L_1^t\in
\d x,\int_t^{\infty}\d b\, L_1^b\in \d z\right)E[G(\Delta T_{[0,x]})|T_x=z].$$
Call $J(x,z)$ the expectation in the integral on the right hand side. Dividing
by $t$ and letting $t\downarrow0$
should yield the generator of the $\R_+^2$-valued
process $((L^t_1,\int_t^{\infty}\d b L_1^b),t\geq 0)$, evaluated at the function
$J$ and at the starting point $(0,1)$. Now, we interpret (see Sect.\ 
\ref{secCSBP} for definitions) the process
$(L_1^t,t\geq 0)$ under $N^{(1)}$ as the $\alpha$-CSBP 
conditioned both to start at
$0$ and stay positive, and 
to have a total progeny equal to $1$. It is thus heuristically a Doob
$h$-transform of the initial CSBP with harmonic function $h(x)=x$, and 
conditioned to come back near $0$ when its integral comes near $1$.  
Now as a consequence of Lamperti's time-change between CSBP's and 
L\'evy processes, the generator of the
CSBP started at $x$ is $x{\cal L}(x,\d y)$ where ${\cal L}$ is the generator of
the stable L\'evy process with index $\alpha$:
$${\cal
L}f(x)=\int_0^{\infty}\frac{C_{\alpha}\d y}{y^{\alpha+1}}
(f(x+y)-f(x)-yf'(x)),$$
where $f$ stands for a generic function in the Schwartz space. 
This, together with well-known properties for generators of $h$-transforms
allows to conjecture that the generator ${\cal L}'$ of the CSBP conditioned to
stay positive and started at $0$ is given by 
$${\cal L}'f(0)=\int_0^{\infty}\frac{C_{\alpha}\d y}{y^{\alpha}}(f(y)-f(0)),$$
for a certain class of nice functions $f$. On the other hand, conditioning to
come back to $0$ when the integral attains $1$ should introduce the term
$q_y(1)$ (recall its definition (\ref{ballotc})) 
in the integral with a certain coefficient, since the total progeny
of a CSBP started at $y$ is equal in law to $T_y$, as a consequence of 
Ray-Knight's theorem. To be more accurate, 
the CSBP starting at $y$ and conditioned to stay positive
should be in $[\eps,\eps+d\eps]$ when its integral equals $1$
with probability close to $\eps y^{-1}q_y(1)d\eps$. 
Indeed, by the conditioned form of Lamperti's theorem of \cite{lambert01} and 
to be recalled below, this
is the same as the probability that the L\'evy process started at $y$ and 
conditioned to stay positive is 
in $[\eps,\eps+d\eps]$ at time $1$. Then,
$$P_y(X_1\in d\eps| T_0>1)= 
\eps y^{-1}P_y(X_1\in d\eps, T_0>1)
\build\sim_{\eps\downarrow0}^{} \eps y^{-1}q_y(1-\eps) d\eps.$$
This, thanks to Lemma \ref{msdisl}, allows to conjecture the form of the
dislocation measure as 
$$\nu_-(G)=C\int_0^{\infty}\frac{\d y\, q_y(1)}{y^{\alpha+1}}E[G(\Delta
T_{[0,y]})|T_y=1]$$ 
for some $C>0$, that can be shown to be equal to $\alpha D_{\alpha}$ with 
some extra care, but we do not need it at this point. It is then easy to reduce
this to the form of Theorem \ref{T1}: 
by using the scaling identities and changing variables $u=y^{-\alpha}$, 
we have that 
\begin{eqnarray*}
\int_0^{\infty}\frac{\d y\, q_y(1)}{y^{\alpha+1}}E[G(\Delta T_{[0,y]})|T_y=1]
&=&\int_0^{\infty}\frac{\d y\, q_1(y^{-\alpha})}{y^{2\alpha-1}}
E[G(y^{\alpha}\Delta T_{[0,1]})
|y^{\alpha}T_1=1]\\
&=&\int_0^{\infty}\alpha^{-1}\d u\, u\, q_1(u) E[G(u^{-1}\Delta T_{[0,1]})|T_1=u]
\\
&=&\alpha^{-1}E[T_1 G(T_1^{-1}\Delta T_{[0,1]})],
\end{eqnarray*}
as wanted.

This very rough 
program of proof could probably be ``upgraded'' to a real rigorous proof, but
the technical difficulties on generators of processes would undoubtedly make
it quite involved. We are going to use a path that uses more the structure of
the stable tree.

%Suppose the starting point $L_1^0$ was
%small but positive instead of $0$. Then $(L^t_1,t\geq 0$ would be interpreted,
%at the light of the discussion of section \ref{secsttree}, as the ``stable'' 
%CSBP starting at $\eta>0$, and conditioned to have a total integral equal to
%$1$. 

\subsection{Erosion and first properties of the dislocation measure}
\label{erod}

From this section on, $F^-$ is exclusively defined under $N^{(1)}$, so that we
may use the nicer notations $P(F^-(t)\in
\d\bs)$ or $E[G(F^-(t))]$ instead of $N^{(1)}(F^-(t)\in \d\bs)$ or 
$N^{(1)}(G(F^-(t)))$ if there is no ambiguity. 

\begin{lmm}\label{erofm}
The erosion coefficient $c$ of $F^-$ is $0$, and the dislocation 
measure $\nu_-(d s)$
charges only $\{s\in S :\ 
\sum_{i=1}^{+\infty}s_i=1\}$.
\end{lmm}

\begin{proof}
We will follow the analysis of Bertoin 
\cite{bertsfrag02}, using 
the law of the time at which a tagged fragment vanishes. 
Let $U$ be uniform on $(0,1)$ and independent of the
height process of the stable tree. Recall the definition of $F^-(t)$ out 
of the open set $I_-(t)$ and let $\lambda(t)=|I^*(t)|$ be the size of the 
interval $I_-^*(t)$ of $I_-(t)$ that contains $U$. As in Sect.\ 
\ref{ssfrag}, if
we define
$$a(t)=\inf \left\{u\geq 0\ :\ \int_0^u\lambda(v)^{
1/\alpha-1}\d v>t\right\}\quad , \quad t\geq 0,$$
then $(-\log(\lambda(a(t))),t\geq 0)$ 
is a subordinator with Laplace exponent
\begin{equation}\label{phi}
\Phi(r)=c(r+1)+\int_S\left( 
1-\sum_{n=1}^{+\infty}s_n^{r+1}\right)
\nu_-(d s).
\end{equation}
Moreover, if $\xi=H_U$ is the lifetime of the tagged fragment, then
\begin{equation}\label{Carmona}
E[\xi^k]=\frac{k!}{\prod_{i=1}^{k}\Phi\left( i\left( 1-\frac{1}{\alpha}
\right)\right)}.
\end{equation}
For the computation we are going to use Proposition \ref{skelstruc}. Recall
that $\chi_s(u)$ is characterized by its Laplace transform
\begin{equation}\label{densitystable}
\int_0^{+\infty}e^{-\mu u}\chi_s(u)\d u=\exp(-s\mu^{1-1/\alpha}).
\end{equation}
We may now compute the moments of $\xi$. We have 
$$E[\xi^k]=\int_0^{+\infty}h^k\alpha\Gamma\left( 1-\frac{1}{\alpha}\right) 
\chi_{\alpha
h}(1)\d h=\frac{\Gamma\left( 1-\frac{1}{\alpha}\right)}{\alpha^k}
\int_0^{+\infty} s^k 
\chi_s(1)\d s.$$ 
To compute this we use (\ref{densitystable}) and Fubini's theorem to get
$$\int_0^{+\infty}\d u e^{-\mu u} \int_0^{+\infty}\d s \chi_s(u)s^k=
\int_0^{+\infty} s^k
\exp(-s\mu^{1-1/\alpha})\d s=\frac{k!}{\mu^{(k+1)
(1-1/\alpha)}},$$
and then the last term above is equal to
$$\frac{k!}{\Gamma\left( (k+1)\left( 1-\frac{1}{\alpha}\right)\right)}
\int_0^{+\infty}
\d u\, e^{-\mu u}
u^{(k+1)(1-1/\alpha)-1}.$$
Inverting Laplace transforms and taking $u=1$ thus give
$$\int_0^{+\infty} s^{k}\chi_s(1)\d s=
\frac{k!}{\Gamma\left((k+1)\left( 1-\frac{1}{\alpha}\right)\right)},$$
hence we finally get
$$E[\xi^k]=\frac{k!\Gamma
\left( 1-\frac{1}{\alpha}\right)}{\alpha^k\Gamma\left((k+1)\left(
1-\frac{1}{\alpha}\right)\right)}.$$
Using (\ref{Carmona}) we now obtain that 
$$\Phi\left( k\left(
1-\frac{1}{\alpha}\right)\right)=
\alpha\frac{\Gamma\left((k+1)\left( 1-\frac{1}{\alpha}
\right)\right)}{\Gamma\left(
k\left(1-\frac{1}{\alpha}\right)\right)},\ \ \ \ k=1,2,\ldots$$
%It is not difficult to write this in L\'evy-Khintchine form, that is,
%$$\alpha\frac{\Gamma\left((k+1)\left( 1-\frac{1}{\alpha}
%\right)\right)}{\Gamma\left(
%k\left(1-\frac{1}{\alpha}\right)\right)}=\int_0^{\infty}
%\frac{\left( 1-\frac{1}{\alpha}\right) e^x \d x}{\Gamma\left(
%1+\frac{1}{\alpha}\right)(e^x-1)^{2-1/\alpha}}\left(1-\exp\left(-xk\left(
%1-\frac{1}{\alpha}\right)\right)\right)$$
%, that
%$\Phi$ is  uniquely determined by the above values, and
Thus, for $r$ of the form $k(1-1/\alpha)$, 
\begin{equation}\label{phir}
\Phi(r)=\alpha\frac{\Gamma\left( r +1-\frac{1}{\alpha}\right)}{\Gamma (r)}
=\frac{r}{\Gamma\left( 1+\frac{1}{\alpha}\right)} B\left(
r+1-\frac{1}{\alpha},\frac{1}{\alpha}\right).
\end{equation}
It is not difficult, using the integral representation of the function $B$,
then changing variables and integrating by parts, to write this in
L\'evy-Khintchine form, that is,
\begin{equation}\label{levkhform}
\frac{r}{\Gamma\left( 1+\frac{1}{\alpha}\right)} B\left(
r+1-\frac{1}{\alpha},\frac{1}{\alpha}\right)=\int_0^{\infty}
\d x\frac{\left( 1-\frac{1}{\alpha}\right) e^x }{\Gamma\left(
1+\frac{1}{\alpha}\right)(e^x-1)^{2-1/\alpha}}\left(1-e^{-xr}\right),
\end{equation}
and it follows that (\ref{phir}) remains true for every $r\geq 0$, hence
generalizing Equation (12) in \cite{bertsfrag02} in the Brownian case. It also
gives the formula 
$$L(\d x)=\frac{\left( 1-\frac{1}{\alpha}\right) e^x \d x}{\Gamma\left(
1+\frac{1}{\alpha}\right)(e^x-1)^{2-1/\alpha}}$$
for the L\'evy measure $L(\d x)$ of $\Phi$, hence generalizing Equation (11) in
\cite{bertsfrag02}.

To conclude, we just notice that $\Phi(0)=0$, which by (\ref{phi}) gives 
both $c=0$ and
$\int_{S}\nu_-(\d\bs)(1-\sum_{i=1}^{\infty}s_i)=0$, implying 
the result. \cq
\end{proof}

%\rem 
%It has to be noticed that the nature of the fragmentation $F^-$ is quite 
%different from that of the ``Brownian fragmentation'' depicted in 
%\cite{bertsfrag02}, since the latter had a ``binary'' fragmentation
%measure whereas the former allows only fragmentations with an infinite number
%of fragments. This could sound paradoxical since our analysis should
%generalize that of \cite{bertsfrag02}, and the results therein should be
%obtained by the substitution $\alpha=2$ in our result. Notice, however, 
%that in this case, $C_-=0$! This means that in the Brownian case $\alpha=2$, 
%we have to be more careful in passing to the limit above. The point is that 
%in this case, the renormalization in front of 
%$\P((\Delta T_{[0,\ell]})\in d{\bf s}|T_{\ell}=1)$ must 
%force the two first components of $(\Delta T_{[0,\ell]})$
%to have a sum close to $1$, a behavior that applies to every stable
%subordinator as explained in \cite{MS03}. 

\subsection{Dislocation measure}

The dislocation measure of $F^-$ will now be obtained by explicitly computing
the law of the first fragmentation of the fragments marked by $n$
independent uniform variables $U_1,\ldots,U_n$ on $(0,1)$, as explained in
Sect.\ \ref{ssfrag}. 
This is going to be a purely combinatorial computation based on
the first formula of Proposition \ref{skelstruc}. What we want to compute is
the law of the partition of $n$ induced by the partition $I_-(t_n)$ and the
variables $U_1,\ldots,U_n$ at the time $t_n$ when they are first separated. In
terms of the stable tree described in section \ref{secsttree}, the probability
$\rho_-(n)(\{\pi_n\})$ 
that the partition induced by $I_-(t_n)$ equals some non-trivial partition
$\pi_n$ of $[n]$ with blocks $A_1,\ldots,A_k$ having sizes $n_1,\ldots,n_k$
with sum $n$ ($n,k\geq 2$) is simply the probability that the
skeleton of the marked tree $\vartheta_n$ is such that the root has
out-degree $k$, and the $k$ trees that are rooted at the children of the root
have $n_1,n_2,\ldots,n_k$ leaves, times $n_1!\ldots n_k!/n!$, which is the 
probability that the labels of these leaves, inherited from the 
sample $(U_i,1\leq i\leq n)$, induce the right partition. 
Let $\bT^*_{n_1,\ldots,n_k}$ be the
set of trees of $\bT^*_n$ that have this last property.
For $x\geq 0$ and $n\geq 0$ we denote by $[x]_n$ the quantity 
$\prod_{i=0}^{n-1}(x+i)=\Gamma(x+n)/\Gamma(x)$. 

\begin{lmm}\label{rhosttree}
Let $\pi_n$ be a partition of $[n]$ with $k\geq 2$ blocks having sizes 
$n_1,n_2,\ldots,n_k$. Then
$$\rho_-(n)(\{\pi_n\})=\frac{D_{\alpha}\Gamma(k-\alpha)}{\alpha^{k}
\Gamma\left(n-\frac{1}{\alpha}\right)}\prod_{i=1}^k\left[
1-\frac{1}{\alpha}\right]_{n_i-1}.$$
\end{lmm}

\begin{proof}
Recall that we want to compute the probability that the skeleton of the 
marked tree $\vartheta_n$ has a root with $k$ children, and the 
fringe subtrees spanned by these children are trees of $\bT_{n_i}^*$
for $1\leq i\leq k$. The fact that the first displayed quantity in 
Proposition \ref{skelstruc} defines a probability on $\bT_n^*$ implies
\begin{eqnarray*}
\sum_{\TT\in\bT_n^*}\prod_{v\in {\cal N}_{\TT}}\frac{|(\alpha-1)
(\alpha-2)\ldots (\alpha-c_v(\TT)+1)|}{c_v(\TT)!}&=&\frac{(\alpha-1)(2\alpha-1)
\ldots((n-1)\alpha-1)}{n!}\\
&=&\frac{\alpha^{n-1}}{n!}\left[1-\frac{1}{\alpha}\right]_{n-1}.
\end{eqnarray*}
Now we compute 
\begin{eqnarray*}
\rho_-(n)(\{\pi_n\})&=&\sum_{\TT\in\bT^*_{n_1,\ldots,n_k}}
\frac{n!n_1!\ldots n_k!}{\alpha^{n-1}\left[1-\frac{1}{\alpha}\right]_{n-1}n!}
\prod_{v\in {\cal N}_{\TT}}\frac{|(\alpha-1)
(\alpha-2)\ldots (\alpha-c_v(\TT)+1)|}{c_v(\TT)!}\\
&=&\frac{n_1!\ldots n_k!
|(\alpha-1)(\alpha-2)\ldots(\alpha-k+1)|}{\alpha^{n-1}k!
\left[1-\frac{1}{\alpha}\right]_{n-1}}\\
&\times& \sum_{\TT\in\bT^*_{n_1,\ldots,n_k}}\prod_{
v\in {\cal N}_{\TT}\setminus\{{\rm root}\}}\frac{|(\alpha-1)
(\alpha-2)\ldots (\alpha-c_v(\TT)+1)|}{c_v(\TT)!}\\
&=&\frac{(\alpha-1)\Gamma(k-\alpha)\Gamma\left(1-\frac{1}{\alpha}
\right)}{k!\alpha^{n-1}\Gamma(2-\alpha)\Gamma
\left(n-\frac{1}{\alpha}\right)}\\
&\times& k!n_1!\ldots n_k! \prod_{i=1}^k \sum_{\TT\in 
\bT_{n_i}^*}\prod_{v\in{\cal N}_{\TT}}\frac{|(\alpha-1)
(\alpha-2)\ldots (\alpha-c_v(\TT)+1)|}{c_v(\TT)!},
\end{eqnarray*}
where the last equality stems from the definition of 
$\bT_{n_1,\ldots,n_k}^*$, and where the factor
$k!$ appears because the $k$ fringe subtrees spanned by the 
sons of the root may appear in any order. By the first 
formula of the proof this now reduces to 
$$\rho_-(n)(\{\pi_n\})=\frac{D_{\alpha}\Gamma(k-\alpha)\prod_{i=1}^k n_i!}
{\alpha^{n}\Gamma
\left(n-\frac{1}{\alpha}\right)}\prod_{i=1}^k\frac{\alpha^{n_i-1}}{n_i!}
\left[1-\frac{1}{\alpha}\right] _{n_i-1}, $$
giving the result.
\cq
\end{proof}

Comparing with Lemma \ref{rhocar} implies, since $c=0$, that the dislocation 
measure $\nu_-$ 
of $F^-$ is thus determined up to a multiplicative constant. Since 
we have a conjectured form $D_{\alpha}\nu_{\alpha}$ 
for the dislocation measure $\nu_-$ of $F^-$, 
we just have to compute the quantity $\kappa_-(\pi)$ for 
$\kappa_-$ the exchangeable measure on ${\cal P}$ with frequencies given by 
the conjectured $\nu_-$. Precisely, we have 

\begin{lmm}\label{numpart}
Let $\pi_n$ be a partition of $[n]$ with $k\geq 2$ blocks and block sizes
$n_1,\ldots,n_k$. Then 
$$\kappa_-^n(\{\pi_n\}):=\kappa_-(\{\pi\in{\cal
P}:\pi|_{[n]}=\pi_n\})=\frac{D_{\alpha}\Gamma(k-
\alpha)}{\alpha^{k-1}\Gamma(n-1)}\prod_{i=1}^k\left[1-\frac{1}{\alpha}
\right]_{n_i-1}$$
\end{lmm}

Before proving this we state from (74) in section 6 of 
\cite{pitmanpk02} (notice that the $\alpha$ there is 
our $1/\alpha$):
\begin{prp}\label{piteppf}
Let $\theta>-1/\alpha$, and let $\mu_{\theta}(\d\bs)$ be the measure on $S$ 
$$\frac{\Gamma(\theta+1)}{\Gamma(\alpha\theta+1)} E\left[T_1^{-\theta};
\frac{\Delta T_{[0,1]}}{T_1}\in \d\bs\right].$$
Then $\mu_{\theta}$ is a probability distribution. Moreover, 
let $\pi_n$ be a partition of $[n]$ with non-void block sizes 
$n_1,\ldots,n_k$. Then the probability that the restriction to $[n]$ of the 
exchangeable partition of ${\cal P}$ with frequencies having law $\mu_{\theta}$
is $\pi_n$ is given by 
$$p_{\theta}(n_1,\ldots,n_k)=\frac{[\alpha\theta+1]_{k-1}}{\alpha^{k-1}
[\theta+1]_{n-1}}\prod_{i=1}^k \left[1-\frac{1}{\alpha}\right]_{n_i-1}$$
\end{prp}

\noindent{\bf Proof of Lemma \ref{numpart}. }
The computation of the $\kappa_-^n$ associated with the conjectured 
dislocation measure $\nu_-$ can go through the same lines as in 
\cite{pitmanpk02}, using the explicit densities for size-biased picks 
among the jumps of the subordinator $T$. However, we use 
the following more direct proof. Write $\nu_{\theta}=
D_{\alpha}(\Gamma(\alpha\theta+1)/
\Gamma(\theta+1))\mu_{\theta}$. Recall from the above the 
notation $\kappa_{\bs}(\d\pi)$ for the law of the exchangeable partition
of $\N$ with ranked asymptotic frequencies given by $\bs$. Define
\begin{equation}\label{kaptheta}
\kappa_{\theta}(\d\pi)=\int_S \nu_{\theta}(\d\bs)\kappa_{\bs}(\d\pi)
=D_{\alpha}E\left[T_1^{-\theta}\kappa_{\Delta T_{[0,1]}/T_1}(\d\pi)\right],
\end{equation}
and for $\pi_n$ a partition of $[n]$ with block sizes $n_1,\ldots,n_k$
write $\kappa_{\theta}^n(\{\pi_n\})=\kappa_{\theta}(\{\pi\in {\cal P}:
\pi|_{[n]}=\pi_n\})=
(\Gamma(\alpha\theta+1)/\Gamma(\theta+1))p_{\theta}(n_1,\ldots,n_k)$. 
Notice that when $n,k\geq 2$ and $\bs\in S$, 
we have $\kappa_{\bs}(\{\pi\in{\cal P}:\pi|_{[n]}=\pi_n\})\leq 
1-s_1$. On the other hand, the 
fact that $\nu_-$ integrates $\bs\mapsto 1-s_1$ is easily generalized 
to $\nu_{\theta}$ for $\theta>-1$. We deduce that the map 
$\theta\to\kappa_{\theta}^n(\{\pi_n\})$ is analytic on 
$\{\theta\in \C:{\rm Re}(\theta)>-1\}$. The same holds for 
\begin{equation}\label{intermcalc}
D_{\alpha}\frac{\Gamma(\alpha\theta+1)}{\Gamma(\theta+1)}
p_{\theta}(n_1,\ldots,n_k)=\frac{D_{\alpha}\Gamma(\alpha\theta+k)}{
\alpha^{k-1}\Gamma(\theta+n)}\prod_{i=1}^k 
\left[1-\frac{1}{\alpha}\right]_{n_i-1}
\end{equation}
provided $k\geq 2$, so the limits as $\theta\in\R \downarrow -1$ 
of $(\ref{kaptheta})$ and of (\ref{intermcalc}) coincide. Using 
a dominated and monotone convergence argument to get the 
$\theta\downarrow -1$ limit in (\ref{kaptheta}), we finally obtain 
$$\kappa_-^n(\{\pi_n\})=\frac{D_{\alpha}\Gamma(k-\alpha)}{
\alpha^{k-1}\Gamma(n-1)}\prod_{i=1}^k 
\left[1-\frac{1}{\alpha}\right]_{n_i-1},$$
as wanted. \cq

\noindent{\bf Remark. }
By analogy with the EPPF (exchangeable partition probability function)
that allows to characterize the law of exchangeable partitions, expressions
such as in Lemma \ref{numpart} could be called ``exchangeable partition
distribution functions'', as they characterize $\sigma$-finite exchangeable
measures on the set of partitions of $\N$. The expression in Lemma
\ref{numpart} should be interpreted as an EPDF for a generalized 
$(1/\alpha,\theta)$ partition (see \cite{pitmancsp02}), for 
$\theta=-1$. One certainly could imagine more general exchangeable partitions
as $\theta$ goes further in the negative axis: this would impose more and 
more stringent constraints on the number of blocks of the partitions.

Therefore, we obtain that
$$\kappa_-^n=\alpha(\Gamma(n-1/\alpha)/\Gamma(n-1))\rho_-(n)$$
on the set of non-trivial partitions of $[n]$. Lemma \ref{rhocar} implies that
the dislocation measure of $F^-$ is equal to the conjectured $\nu_-$ up to a
multiplicative constant. We are going to recover the missing information with
the help of the computation of $\Phi$ above. 

\subsection{The missing constant}

In this section, we compute the Laplace exponent
$\Phi$ of the subordinator $-\log(\lambda(a(\cdot)))$ 
of Sect.\ \ref{erod}, whose value is indicated in (\ref{phir}),
directly from formulas (\ref{phi}) and (\ref{nupm}). 
Let
$$\Phi_0(r)=\int_S\left(1-\sum_{n=1}^{\infty}s_n^{r+1}\right)\nu_{-}(\d\bs),$$
where $\nu_-$ is the measure given in Theorem \ref{T1}.
If we can prove that $\Phi_0(r)=\Phi(r)$ for every $r\geq 0$, we will
therefore have established that the normalization of $\nu_-$ is the
appropriate one. By (\ref{nupm}), 
\begin{eqnarray*}
\Phi_0(r)&=&D_{\alpha}E\left[T_1\left(1-\sum_{0\leq x\leq 1}
\left(
\frac{\Delta T_x}{T_1}\right)^{r+1}\right)\right]\\
&=&D_{\alpha}\int_0^{\infty}\d u\,u\, q_1(u)E\left[1-\sum_{
0\leq x\leq 1}\left(\frac{
\Delta T_x}{u}\right)^{r+1}\bigg|T_1=u \right]\\
&=&D_{\alpha}\int_0^{\infty}\d u\, u \, q_1(u)E\left[1-\left(
\frac{\Delta^*_1}{u}\right)^r
\right]
\end{eqnarray*}
where $\Delta^*_1$ is a size-biased pick from the jumps of $T_x$, for 
$0\leq x\leq 1$, conditionally on $T_{1}=u$. Using formula 
(\ref{sizeb}) and recalling that $T$ has L\'evy measure 
$c_{\alpha}x^{-1-1/\alpha}\d x$, we can write
\begin{eqnarray*}
\Phi_0(r)&=&D_{\alpha}\int_0^{\infty}\d u\, u\,q_1(u)\int_0^u 
\d x(1-(x/u)^r)
\frac{c_{\alpha} q_1(u-x)}{ux^{1/\alpha}q_1(u)}\\
&=&D_{\alpha}\int_0^{\infty}\d u \int_0^1 \d y\, c_{\alpha}u^{1-
1/\alpha} q_1(u(1-y))
\frac{1-y^r}{y^{1/\alpha}}\\
&=&D_{\alpha}\int_0^1 \d y \frac{c_{\alpha}(1-y^r)}{y^{1/
\alpha}(1-y)^{2-1/\alpha}}
\int_0^{\infty}\d u\, u^{1-1/\alpha}q_1(u)
\end{eqnarray*}
as obtained by Fubini's theorem, and linear changes of variables. 
The integral in $\d u$
equals $\E[T_1^{1-1/\alpha}]$, which is 
$\Gamma(2-\alpha)/\Gamma(1/\alpha)$ by standard results using Laplace
transform. 
Using the expressions for $D_{\alpha}$, $c_{\alpha}$ and the identity
$\alpha^{-1}\Gamma(1/\alpha)=\Gamma(1+1/\alpha)$, it remains to compute the
quantity
$$\frac{1-\frac{1}{\alpha}}{\Gamma\left(1+\frac{1}{\alpha}\right)}\int_0^1\frac{\d
y\,  y^{-1/\alpha}(1-y^r)}{(1-y)^{2-1/\alpha}}.$$
But this is exactly the expression
(\ref{levkhform}) after changing variables $y=e^{-x}$, and it is thus equal to
$r B(r+1-1/\alpha,1/\alpha)/\Gamma(1+1/\alpha)$, which is (\ref{phir}) as
wanted, thus completing the proof of Theorem \ref{T1}.

\section{Small-time asymptotics}\label{asymzero}

In this section we study the asymptotic behavior of $F^-$ for small times. 
Precisely, let $M(t)=\sum_{i\geq 1}F^-_i(t)$ 
denote the total mass of $F^-$ at time $t$. Let $(Y_x,x\geq 0)$ denote 
an $\alpha$-CSBP, started at 
$0$ and conditioned to stay positive. See the following section for the 
definitions. We have the following 
result, that generalizes and mimics somehow results from 
\cite{jpda98sac,berest02,MS03}. However, these results dealt with self-similar 
fragmentations with positive indices, and also, the occurrence of the 
randomization introduced by $Y_1$ below is somehow unusual.

\begin{prp}\label{smasymfm}
The following convergence in law holds:
$$t^{\alpha/(1-\alpha)}(M(t)-F^-_1(t),F^-_2(t),F^-_3(t),\ldots)\build\to_{t
\downarrow 0}^{d} (T_{Y_1},\Delta_1,\Delta_2,\ldots)$$
where $T$ is the stable $1/\alpha$ subordinator as above, independent of $Y$, 
and $\Delta_1,\Delta_2,\ldots$ are the jumps of $(T_x,0\leq x\leq Y_1)$ ranked
in decreasing order of magnitude. 
\end{prp} 

For this we are going to use the following lemma, which resembles the result
of Jeulin in \cite{jeulin80} relating a scaled normalized Brownian excursion
and a 3-dimensional Bessel process. The proof is  postponed to the
following section. Recall that
$(L_1^t,t\geq 0)$ stands for the local time of the height process up to time
$1$. 

\begin{lmm}\label{asymind}
The following convergence in law holds:
$$\mbox{Under }N^{(1)},\quad 
(t^{1/(1-\alpha)}L_1^{tx},x\geq0)
%t^{\alpha/(1-\alpha)}\int_0^t \d s\,L_1^s\right)
\build\to_{t\downarrow0}^{d}(Y_x,x\geq 0),$$
%, \int_0^1 \d s\, Y_s\right),$$
and this last limit is independent of the initial process $(L_1^t,t\geq 0)$.
In particular, $t^{1/(1-\alpha)}L_1^t$
converges in distribution to
$Y_1$ as $t\downarrow0$.
\end{lmm}

In the sequel let $(y_t,\ov{y}_t)$ have the law of 
$(L_1^t,\int_t^{\infty} \d b L_1^b)$ under $N^{(1)}$.

\noindent{\bf Proof of Proposition \ref{smasymfm}. }
Following the method of Aldous and Pitman \cite{jpda98sac}, we are actually
going to prove that for every $k$, 
\begin{equation}\label{recureq}
t^{\alpha/(1-\alpha)}(M(t)-F^*_1(t),F^*_2(t),F^*_3(t),\ldots,F^*_k(t))
\build\to_{t\downarrow 0}^{d}
(T_{Y_1},\Delta^*_1,\Delta^*_2,\ldots,\Delta^*_{k-1}),
\end{equation}
for every $k\geq 1$, 
where the quantities with the stars are the size-biased quantities associated
with the ones of the statement, and this is sufficient. We are going to proceed
by induction on $k$. To start the induction, let $g$ be a continuous
function with compact support
and write, using Lemma \ref{permanform}, Proposition \ref{semigrmoins}, 
then changing variables and using scaling identities,
\begin{eqnarray}
E[g(t^{\alpha/(1-\alpha)}(M(t)-F^*_1(t)))]&=&E\left(\int_0^{\ov{y}_t}
\d u\frac{c_{\alpha}\, y_t\, q_{y_t}(\ov{y}_t-u)}{
\ov{y}_t\, 
u^{1/\alpha}\, q_{y_t}(\ov{y}_t)} g(t^{\alpha/(1-\alpha)}(\ov{y}_t-u))
\right)\nonumber\\
&=&E\left(\int_0^{t^{\alpha/(1-\alpha)}\ov{y}_t}\!\!
\d v\frac{t^{\alpha/(\alpha-1)
}\, 
c_{\alpha}\, y_t\,
q_1\left(\frac{v}{t^{\alpha/(1-\alpha)}y_t^{\alpha}}\right)}{
(\ov{y}_t-t^{\alpha/(\alpha-1)}v)^{1/\alpha }\, \ov{y}_t\,
q_1\left(\frac{\ov{y}_t}{y_t^{\alpha}}\right)}g(v)\right).\nonumber
\end{eqnarray}
By making use of  Skorokhod's representation theorem, we may suppose that the 
convergence of $(t^{1/(1-\alpha)}y_t,t^{\alpha/(1-\alpha)}\ov{y}_t)$ to 
$(Y_1,\infty)$ is almost-sure. Now
the integral inside the expectation is the integral according to a probability 
law, hence it is dominated by the supremum of $|g|$, so it suffices 
to show that the integral converges a.s.\ to apply dominated convergence. 
For almost every $\omega$, there exists $\eps$ such that if $t<\eps$, 
$t^{\alpha/(1-\alpha)}\ov{y}_t(\omega)>K$ where $K$ is the right-end of the 
support of $g$. For such an $\omega$ and $t$, the integral is thus 
\begin{eqnarray*}
\lefteqn{\int_0^K \d v\, g(v)\frac{c_{\alpha}t^{\alpha/(\alpha-1)}y_t 
q_1(v(t^{1/(1-\alpha)}y_t)^{-\alpha})}{\ov{y}_t^{1+1/\alpha}(1-t^{\alpha/
(\alpha-1)}v/\ov{y}_t)^{1/\alpha}q_1(\ov{y}_t y_t^{-\alpha})}}\\
&\leq&
M\frac{t^{\alpha/(\alpha-1)}y_t}{\ov{y}_t^{1+1/\alpha}q_1(\ov{y}_t y_t^{
-\alpha})}\int_0^K \d v\, q_1\left(\frac{v}{ t^{\alpha/(1-\alpha)
}y_t^{\alpha}}\right)
\end{eqnarray*}
for some constant $M$ not depending on $t$. Now we use the fact
%By reordering slightly this sum, we are thus taking the expectation of the 
%function 
%\begin{equation}\label{funcG}
%G(x,y)=\int_0^{\infty}\d v g(v)\frac{c_{\alpha}\ind_{[0,y]}(v)x
%q_1(v/x^{\alpha})}{ y^{1+1/\alpha}(1-v/y)^{1/\alpha}q_1(y/x^{\alpha})}=
%\int_0^{1}\d v g(vy)\frac{c_{\alpha}\, x\, q_1(vy/x^{\alpha})}{
%y^{1/\alpha}(1-v)^{1/\alpha}q_1(y/x^{\alpha})}
%\end{equation}
%under the law of the pair
%$(t^{1/(1-\alpha)}y_t,t^{\alpha/(1-\alpha)}\ov{y}_t)$,  which converges is law
%to $(Y_1,\infty)$. Thus if one can show that $G$ is continuous bounded, we
%will be able to pass to the limit. To justify both  we write
%$$g(vy)\frac{x q_1(vy/x^{\alpha})}{
%y^{1/\alpha}(1-v)^{1/\alpha}q_1(y/x^{\alpha})}\leq \| g\|_{\infty}
%\frac{q_1(v\sigma)}{\sigma^{1/\alpha}(1-v)^{1/\alpha}q_1(\sigma)}$$
%where $\sigma=y/x^{\alpha}$. It is then easy by dominated convergence to see
%that $G$ is continuous on $(0,\infty)^2$. It thus remains to see that it has a
%finite $\limsup$ as $\sigma\to0$ and $\sigma\to\infty$. For this we use the
from \cite{skorokhod} that $q_1$ is bounded and
%there exists constants $A,B>0$ such that 
%$$q_1(x)\build\sim_{x\to0+}^{} A x^{-1-1/(2(\alpha-1))}\exp(-B
%x^{-1/(\alpha-1)})\quad
%\mbox{ and }\quad
$$q_1(x)\build=_{x\to\infty}^{}c_{\alpha}x^{-1-1/\alpha}+O(x^{-1-2/\alpha}).$$
%This implies that for large $\sigma$, there exist constants $M,M',M''$ such
%that 
%\begin{eqnarray*}
%G(x,y)&\leq& \int_0^1 \d v\left(\frac{\ind_{\sigma v\leq
%1}M}{\sigma^{1/\alpha}q_1(\sigma)}+\frac{\ind_{\sigma v\geq 1}
%M'}{\sigma^{1/\alpha}v^{1+1/\alpha}(1-v)^{1/\alpha}}\right)\\
%&\leq
%&\int_0^{1/\sigma}\d v M''\sigma+\int_{1/\sigma}^1\frac{M'}{
%\sigma^{1/\alpha}v^{1+1/\alpha}(1-v)^{1/\alpha}}
%\end{eqnarray*}
%The first integral is thus bounded by $M''$, and it is easy to see that the
%behavior of the second at the bound $1/\sigma$ is as
%$\sigma^{1/\alpha}/\sigma^{1/\alpha}$. So the two are bounded when
%$\sigma\to\infty$. When $\sigma$ is small, there exists $M$
%such that 
%\begin{eqnarray*}
%G(x,y)&\leq&
%\frac{M}{\sigma^{1/\alpha}}\int_0^1\d v\frac{
%\exp(-B\sigma^{-1/(\alpha-1)}(v^{-1/(\alpha-1)}-1))}{v^{1+1/2(\alpha-1)}
%(1-v)^{1/\alpha}}\\
%&=&M\sigma^{1+1/2(\alpha-1)}\exp(B\sigma^{1/(1-\alpha)})
%\int_0^{\sigma}\d v\frac{\exp(-Bv^{1/(1-\alpha)})}{v^{1+1/2(\alpha-1)}(\sigma-v
%)^{1/\alpha}}.
%\end{eqnarray*}
%from this it is easy to show that the whole goes to $0$ as $\sigma\to0$, by
%splitting the integral according to $0\leq v\leq \sigma/2$ and $\sigma/2\leq
%v\leq \sigma$. This by (\ref{expg},\ref{funcG}) justifies the limit 
%$$E[g(t^{\alpha/(1-\alpha)}(M(t)-F^*_1(t)))]\to E\left[\int_0^{\infty}\d v
%g(v)\frac{q_1(v/Y_1^{\alpha})}{Y_1^{\alpha}}\right]=E[g(T_{Y_1})],$$
%where we have again used the asymptotic behavior of $q_1(x)$ as $x\to\infty$.
%and the scaling identity. 
This allows to conclude by dominated convergence that the 
integral a.s.\ goes to 
$$\int_0^K \d v\, g(v)\frac{q_1(v/Y_1^{\alpha})}{Y_1^{\alpha}}=
\int_0^K \d v\, g(v) q_{Y_1}(v), $$
and by dominated convergence its expectation converges to the expectation of 
the above limit, that is $E[g(T_{Y_1})]$. 

To implement the recursive argument, suppose that
(\ref{recureq}) holds for some $k\geq 1$. Let $g$ and $h$ be continuous
bounded functions on $\R_+$ and $\R_+^k$ respectively. Denote by
$(y_t,\ov{y}_t,\Delta_1(t),\Delta_2(t)\ldots)$ a sequence with the
same law as $(L_1^t,\int_t^{\infty}\d s L_1^s,\Delta T'_{[0,L_1^t]})$ given
$T'_{L_1^t}=
\int_t^{\infty}\d s L_1^s$, where $L_1$ is taken under $N^{(1)}$ and
$T'$ is a stable $1/\alpha$ subordinator, taken
independent of $L$. Last, let $\Delta^*_1(t),\Delta^*_2(t),\ldots$ 
be the size-biased permutation associated with 
$\Delta_1(t),\Delta_2(t),\ldots$.
By Proposition
\ref{semigrmoins}, conditioning and using Lemma \ref{permanform} we have 
\begin{eqnarray*}
\lefteqn{E[g(t^{\alpha/(1-\alpha)}
F^*_{k+1}(t))h(t^{\alpha/(1-\alpha)}(M(t)-F^*_1(t),F^*_2(t),\ldots,
F^*_k(t)))]}\\
&=&E\bigg[h(t^{\alpha/(1-\alpha)}(\ov{y}_t-\Delta^*_1(t),\Delta^*_2(t),\ldots,
\Delta^*_k(t)))\\
&\times&\int_0^{\ov{y}_t-\sum_{i=1}^k\Delta^*_i(t)}\d u
\, g(t^{\alpha/(1-\alpha)}u)
\frac{
c_{\alpha}\, y_t\, q_{y_t}\left(\ov{y}_t-\sum_{i=1}^k\Delta^*_i(t)-u\right)}{
u^{1/\alpha}\left(\ov{y}_t-\sum_{i=1}^k\Delta^*_i(t)\right)q_{y_t}\left(
\ov{y}_t-\sum_{i=1}^k\Delta^*_i(t)\right)}\bigg]
\end{eqnarray*}
Similarly as above, we show by changing variables and then using the scaling
identities and the asymptotic behavior of $q_1$ that this converges to
$$E\left[h(T_{Y_1},\Delta^*_1,\ldots,\Delta^*_{k-1})\int_0^{T_{Y_1}-
\sum_{i=1}^k
\Delta^*_i}\!\!\!\d u\, g(u)\frac{c_{\alpha}\,
Y_1\, q_{Y_1}\left(T_{Y_1}-\sum_{i=1}^{k-1}\Delta^*_i-v\right)}{
u^{1/\alpha}(T_{Y_1}-\sum_{i=1}^{k-1}\Delta^*_i)q_{Y_1}\left(T_{Y_1}-
\sum_{i=1}^{k-1}\Delta^*_i\right)}\right]$$
and by Lemma \ref{permanform} 
this is $E[h(T_{Y_1},\Delta^*_1,\ldots,\Delta^*_{k-1})
g(\Delta^*_k)]$. This finishes the proof.
\cq

The same method as that used in this proof can be used to show also that the
rescaled remaining mass $t^{\alpha/(1-\alpha)}(1-M(t))$ converges in
distribution to $\int_0^1Y_v\, \d v$ jointly with the
vector of the proposition.

\section{Some results on continuous-state branching processes}\label{secCSBP}

In this section we develop the material needed to prove Lemma \ref{asymind}.
In  the course, we will give an analog of Jeulin's theorem \cite{jeulin85}
linking the local time process of a Brownian excursion to another time-changed 
Brownian excursion. To stay in the line of the present paper, we will 
suppose that the laws we consider are associated 
to stable processes, but all of the results (except the proof of Lemma 
\ref{asymind} which strongly uses scaling) can be extended to more general 
L\'evy processes and their associated CSBP's. To avoid confusions, we will
denote by $(Z_t,t\geq 0)$ the different CSBP's we will consider, or to be more
precise, we let $(Z_t,t\geq 0)$ instead of $(X_s,s\geq 0)$ be the canonical
process on $\D([0,\infty))$ when dealing with the laws
$\P_x,\P^{\uparrow}_x,\ldots$ associated to CSBP's. 

\begin{defn}
For any $x>0$, let $\P_x$ be the unique law on $\D([0,\infty))$ 
that makes the canonical process $(Z_t,t\geq0)$ a right-continuous Markov
process starting at $x$ with transition probabilities characterized by 
$$\E[\exp(-\lambda Z_{t+r})|Z_t=y]=\exp(-yu_r(\lambda)),$$
where $u_r(\lambda)=(\lambda^{1-\alpha}+(\alpha-1)r)^{1/(1-\alpha)}$ 
is determined by the equation
$$\int_{u_r(\lambda)}^{\lambda}\frac{\d v}{v^{\alpha}}=r.$$
Then $\P_x$ is called the law of of the $\alpha$-CSBP started at $x$. 
\end{defn}

\noindent{\bf Remark. }
For more general branching mechanisms, the definition of $u_r(\lambda)$ is 
modified by replacing $v^{\alpha}$ by $\psi(v)$, where $\psi$ is the Laplace
exponent of a spectrally positive L\'evy process with infinite variation that 
oscillates or drifts to $-\infty$. 

Recall the setting of section \ref{levyproc}, and let
$P_x$ be law under which $X$ is the spectrally positive 
stable process with Laplace exponent $\lambda^{\alpha}$ 
and started at $x>0$, that is, the law of $x+X$ under $P$. 
Let $E_x$ be the corresponding expectation. 
Define the time-change $(\tau_t,t\geq 0)$ by 
$$\tau_t=\inf\left\{u\geq 0:\int_0^u\frac{\d v}{X_{v\wedge h_0}}>t\right\},$$ 
where $h_0=\inf\{s>0: X_s=0\}$ is the first hitting time of $0$.
This definition makes sense either under the law $P_x$, for $x>0$, or the 
$\sigma$-finite excursion measure $N$ (we will see below that under $N$, 
$\tau$ is not the trivial process identical to $0$).

\begin{thm}\label{jeulinth}
We have the following identities in law: for every $x>0$, 
$$(L_{T_x}^t,t\geq 0) \mbox{ under }P
\build=_{}^{d}(X_{\tau_t},t\geq 0) \mbox{ under }
P_x,$$
and both have law $\P_x$. Moreover, 
$$(L_{\zeta}^t,t\geq 0) \mbox{ under }N \build=_{}^{d}
(X_{\tau_t},t\geq 0)\mbox{ under }N.$$ 
\end{thm}

The first part is already known and is a conjunction of Lamperti's theorem 
and the Ray-Knight theorem mentioned in Sect.\ \ref{secsttree}. We will use it 
to prove the second part. 
First we introduce some notations, which were already 
used in a heuristic way above. 

For $x>0$ one can define the law $P^{\uparrow}_x$ of the stable process started
at $x$ and conditioned to stay positive by means of Doob's theory of 
harmonic $h$-transforms. It is characterized by the 
property 
$$E^{\uparrow}_x[F(X_s,0\leq s\leq K)]=E_x\left[\frac{X_K}{x}F(X_s,0\leq s\leq
K),K<T_0\right]$$
for any positive measurable functional $F$. Here $T_0$ denotes as above the 
first hitting time of $0$ by $X$. 
It can be shown (see e.g.\ \cite{chaumont97}) that $P^{\uparrow}_x$ has a 
weak limit as $x\to 0$, which we call $P^{\uparrow}$, the law of the stable 
process conditioned to stay positive. 

Similarly, we define the CSBP conditioned to stay positive
according to \cite{lambert01}, by letting
$\P_x$ be the law of the CSBP started at $x>0$, then setting 
$$\E^{\uparrow}_x[F(Z_t,0\leq t\leq K)]=\E_x\left[
\frac{Z_K}{x}F(Z_s,0\leq s\leq K)\right].$$
We want to show that a $x\downarrow0$ limit also exists in this case. This
is made possible by the interpretation of \cite{lambert01} of the law
$\P^{\uparrow}_x$ in terms of a CSBP with {\em immigration}. To be concise, we
have 
\begin{lmm}\label{immigr}
For $x>0$, the law $\P_x^{\uparrow}$ is the law of the $\alpha$-CSBP with
immigration function $\alpha\lambda^{\alpha-1}$ and started at $x$. That is,
under $\P^{\uparrow}_x$, $(Z_t,t\geq 0)$ is a Markov process starting at $x$
and with transition probabilities
$$\E^{\uparrow}_x[\exp(-\lambda Z_{t+r})|Z_t=y]=\exp\left(-y
u_r(\lambda)-\int_0^r\alpha u_v(\lambda)^{\alpha-1}\d v\right).$$
\end{lmm}

As a consequence, the laws $\P^{\uparrow}_x$ converge weakly as $x\downarrow0$
to a law $\P^{\uparrow}_0=\P^{\uparrow}$, which is the law of a Markov process
with same transition probabilities and whose entrance law is given by the above
formula, taking $t=y=x=0$. It is also easy that the law $\P^{\uparrow}$ is that
of a Feller process according to the definition for $u_r(\lambda)$. 

It is shown in \cite{lambert01} that Lamperti's correspondence is still 
valid between conditioned processes started at $x>0$:
the process $(X_{\tau_{t}},t\geq 0)$ under the law $P^{\uparrow}_x$ 
has law $\P^{\uparrow}_x$. To be more accurate, the exact statement is that 
if the process $(Z_t,t\geq 0)$ has law $\P^{\uparrow}_x$, 
then the process $(Z_{C_s},s\geq 0)$ has law $P^{\uparrow}_x$ where 
$$C_s=\inf\left\{u\geq 0:\int_0^u\d v Z_v>s\right\},$$ 
but this is the second part of Lamperti's transformation, which is easily 
inverted (see also the comment at the end of the section). 
We generalize this to 
\begin{lmm}\label{amauryth}
The process $(X_{\tau_t},t\geq 0)$ under the law $P^{\uparrow}$ has 
law $\P^{\uparrow}$. 
\end{lmm}
Part of this lemma is that $\tau_t>0$ for every $t$. 

\begin{proof}
For fixed $\eta>0$, let
$$\tau_t^{\eta}=\inf\left\{u:\int_{\eta}^{u\vee\eta}\frac{\d
v}{X_v}>t\right\}.$$ 
This is well defined under $P^{\uparrow}$ since $X_t>0$
for all 
$t>0$ a.s.\ under this 
law.
Then since $\int_{\eta}^{u\vee\eta}\d v/X_v=\int_0^{u-\eta}\d v/X_{\eta+v}$, 
we have that 
$$\tau_t^{\eta}=\eta+\inf\left\{u\geq 0:\int_0^u\frac{\d v}{X_{\eta+v}}>t 
\right\}.$$
That is, $\tau^{\eta}-\eta$ equals
the time-change $\tau$ defined above, but associated to the process 
$(X_{\eta+t},t\geq 0)$ (notice that $h_0$ plays no role here since we are
dealing with processes that are strictly positive on$(0,\infty)$).  Under
$P^{\uparrow}$, this process is  independent of $(X_s,0\leq s\leq \eta)$
conditionally on 
$X_{\eta}$ and has law $P^{\uparrow}_{X_{\eta}}$. 
Hence, 
by Lamperti's identity, conditionally on
$(X_s,0\leq s\leq \eta)$ under $P^{\uparrow}$, 
the process $(X_{\tau_t^{\eta}},t\geq 0)$ has law $\P^{\uparrow}_{X_{\eta}}$.
Hence, for any continuous bounded functional $G$ on the paths defined 
on $[0,K]$ for some $K>0$, 
$$E^{\uparrow}[G(X_{\tau^{\eta}_t}, 0\leq t\leq K)]=E^{\uparrow}[
\E^{\uparrow}_{X_{\eta}}[G(Z_t,0\leq t\leq K)]].$$ 
Now, it is not difficult to see that $\tau^{\eta}$ decreases to the 
limit $\tau$ uniformly on compact sets. Thus, using the right-continuity of 
$X$ on the one hand, and the Feller property on the other (in fact, less than
the Feller property is needed here), we obtain by letting
$\eta\downarrow 0$ in the above identity
$$E^{\uparrow}[G(X_{\tau_t},0\leq t\leq K)]=\E^{\uparrow}[
G(Z_t,0\leq t\leq K)],$$ 
which is the desired identity. In particular, $\tau$ cannot
be identically $0$. \cq
\end{proof}

\noindent{\bf Remark. }
Notice that the fact that the time-change $\tau_t$ is still 
well-defined under the law $P^{\uparrow}$ can be double-checked 
by a law of the iterated logarithm for the law $P^{\uparrow}$. See also the 
end of the section. 

Motivated by the definition in Pitman-Yor \cite{pityor82} for the
excursion measure away from $0$ of continuous diffusions for which 
$0$ is an exit point 
(and initially by It\^o's description of the Brownian excursion
measure linking the three-dimensional Bessel process semigroup to the 
entrance law of Brownian excursions), we now state the following 

\begin{prp}\label{entrlaw}
The process $(L_{\zeta}^t,t\geq 0)$ under the measure $N$ is governed by the 
excursion measure of the CSBP with characteristic $\lambda^{\alpha}$. That 
is, its entrance law $N(L_{\zeta}^t\in \d y)$ for $t>0$ is equal to 
$y^{-1}\P^{\uparrow}(Z_t\in \d y)$ for $y>0$ (and it puts mass $\infty$ on 
$\{0\}$), and given $(L_{\zeta}^u,0\leq u\leq t)$, the process
$(L_{\zeta}^{t+t'},t'\geq 0)$ has law $\P_{L_{\zeta}^t}$. 
\end{prp}

The use of the height process and its local time under $N$, and hence of an
``excursion measure'' associated to the genealogy of CSBP's, snakes and
superprocesses, is a very natural tool, however it does not seem that the above
proposition, which states that this notion of ``excursion measure'' is the
most natural one, has been checked somewhere. However, as noticed 
in \cite{pityor82}, since the point 
$0$ is not an entrance point for the initial CSBP, one cannot define a 
reentering diffusion by sticking the atoms of a Poisson measure with intensity 
given by this excursion measure, because the 
durations are almost never summable. 

\begin{proof}
The law $\P^{\uparrow}(Z_t\in \d y)$ is the weak limit of $\P^{\uparrow}_x(
Z_t\in \d y)=x^{-1}y \P(Z_t\in \d y)$ as $x\to 0$. Since by the properties of
the  CSBP mentioned in section \ref{secsttree}, we have $\E_x[\exp(-\lambda
Z_t)] =\exp(-x u_t(\lambda))$, we obtain 
$$\int_0^{\infty}\frac{\P^{\uparrow}_x(Z_t\in \d y)}{y}(1-e^{-\lambda y})
=\int_0^{\infty}\frac{\P_x(Z_t\in \d y)}{x}(1-e^{-\lambda y})=\frac{1-e^{-x u_t
(\lambda)}}{x}.$$
This converges to $u_t(\lambda)$ as $x\to 0$, and thanks to the proof 
of \cite[Theorem 1.4.1]{duqleg02}, this equals 
$N(1-\exp(-\lambda L_{\zeta}^t))$. This gives the identity of the entrance 
laws. For the Markov property we use excursion theory and 
Ray-Knight's theorem. Let 
$0<t_1<\ldots<t_n<t$, then Markov's property for $(L_{T_1}^t,t\geq 0)$ 
entails that for every $\lambda_1,\ldots,\lambda_n,\lambda\geq 0$,
$$E[\exp(-\sum_{i=1}^n\lambda_iL_{T_1}^{t_i}-\lambda L_{T_1}^t)]
=E[\exp(-\sum_{i=1}^{n-1}\lambda_iL_{T_1}^{t_i}-(\lambda_n+u_{t-t_n}(\lambda)
)L_{T_1}^{t_n})].$$
On the other hand, we may write $L_{T_1}^t=\sum_{0<s\leq 1}
(L_{T_s}^t-L_{T_s-}^t)$ so that the Laplace exponent identity for Poisson
point processes applied to both 
sides of the above displayed expression gives after taking logarithms:
$$N\left(1-\exp\left(-\sum_{i=1}^n\lambda_iL_{\zeta}^{t_i}-
\lambda L_{\zeta}^t\right)\!
\right)
=N\left(1-\exp\left(-\sum_{i=1}^{n-1}\lambda_i
L_{\zeta}^{t_i}-(\lambda_n+u_{t-t_n}(\lambda))L_{\zeta}^{t_n}\right)\!\right),
$$
so that a substraction gives 
\begin{eqnarray*}
\lefteqn{N\left(\exp\left(-\sum_{i=1}^n\lambda_iL_{\zeta}^{t_i}\right)
(1-\exp(-\lambda L_{\zeta}^t))\right)}\\
&=&
N\left(\exp\left(-\sum_{i=1}^n\lambda_iL_{\zeta}^{t_i}\right)
(1-\exp(-u_{t-t_n}(\lambda)L_{\zeta}^{t_n}))\right)\\
&=&N\left(\exp\left(-\sum_{i=1}^n\lambda_iL_{\zeta}^{t_i}\right)
\E_{L_{\zeta}^{t_n}}[1-\exp(-\lambda Z_{t-t_n})]\right). 
\end{eqnarray*}
Hence the Markov property. \cq
\end{proof}

\noindent{\bf Proof of Theorem \ref{jeulinth}. }
It just remains to prove the second statement. For this we let $\eta>0$ and 
we define as above the time change $\tau_{t}^{\eta}$.
Using the Markov property under the measure $N$, we again have that under $N$, 
$(X_{\eta+s},s\geq 0)$ is independent of $(X_s,0\leq s\leq \eta)$ 
conditionally on $X_{\eta}$ and has the law $P_{X_{\eta}}^{h_0}$ of the 
stable process started at $X_{\eta}$ and killed at time $h_0$. Hence, 
by Lamperti's identity, under $N$ and conditionally on
$(X_s,0\leq s\leq \eta)$, 
the process $(X_{\tau_t^{\eta}},t\geq 0)$ has law $\P_{X_{\eta}}$. 
Thus if $\eta<t_1<\ldots<t_n<t$ and if $g_1,\ldots,g_n,g$ 
are positive continuous functions with compact support
that does not contain $0$, 
then
\begin{eqnarray*}
N\left(\prod_{i=1}^n g_i(X_{\tau_{t_i}^{\eta}}) \;
g(X_{\tau_t^{\eta}})\right)&=&\int_0^{\infty} 
N(X_{\eta}\in \d x)\E_x\left[\prod_{i=1}^n 
g_i(Z_{t_i-\eta})\; g(Z_{t-\eta})\right]\\
&=& \int_0^{\infty} N(X_{\eta}\in \d x)\E_x\left[\prod_{i=1}^n 
g_i(Z_{t_i-\eta})\; \E_{Z_{t_n-\eta}}[g(Z_{t-t_n})]\right].
\end{eqnarray*}
As for the CSBP, the entrance law $N(X_{\eta}\in \d x)$ equals 
$x^{-1} P^{\uparrow}(X_{\eta}\in \d x)$ 
for $x>0$. So we recast the last expression 
as
\begin{eqnarray*}
\lefteqn{\int_0^{\infty} 
P^{\uparrow}(X_{\eta}\in \d x)\E_x\left[\frac{\prod_{i=1}^n 
g_i(Z_{t_i-\eta})}{x}\E_{Z_{t_n-\eta}}[g(Z_{t-t_n})]\right]}\\
&=&\int_0^{\infty} 
P^{\uparrow}(X_{\eta}\in \d x)\E^{\uparrow}_x\left[\frac{\prod_{i=1}^n 
g_i(Z_{t_i-\eta})}{Z_{t_n-\eta}}\E_{Z_{t_n-\eta}}[g(Z_{t-t_n})]\right].
\end{eqnarray*}
Now we let $\eta\downarrow 0$, using the right continuity and the Feller 
property of the CSBP, to obtain
$$N\left(\prod_{i=1}^n g_i(X_{\tau_{t_i}}) \;
g(X_{\tau_t})\right)=\E^{\uparrow}\left[\frac{\prod_{i=1}^n 
g_i(Z_{t_i})}{Z_{t_n}}\E_{Z_{t_n}}[g(Z_{t-t_n})]\right].$$
Hence, thanks to Proposition \ref{entrlaw} we obtain that under $N$
the process $(X_{\tau_t},t\geq 0)$ has the same entrance law and Markov 
property as $(L_{\zeta}^t,t\geq 0)$, hence the same law. \cq

\noindent{\bf Proof of Lemma \ref{asymind}. }
Let $G$ be a continuous bounded functional on the paths with
lifetime $K$. We want to show that $N^{(1)}[G(t^{1/(1-\alpha)}L_1^{tx},0\leq
x\leq K)]$ goes to $E^{\uparrow}[G(X_{\tau_x},0\leq x\leq K)]$.
By Theorem \ref{jeulinth}, the process $(L_v^{x},x\geq 0)$ under
$N^{(v)}$ is equal to the process $(X_{\tau_x},x\geq 0)$ under the law 
$N^{(v)}$ for almost every $v$, and we can take $v=1$ by the usual scaling
argument. By \cite{chaumont97}, 
the law $N^{(1)}$ can be obtained as the bridge with
length $1$ of the stable process conditioned to stay positive, and there
exists a positive measurable harmonic 
function $h$ such that for every functional $J$
and every $r<1$, 
$$N^{(1)}[J(X_s,0\leq s\leq r)]=E^{\uparrow}[h(X_r)J(X_s,0\leq s\leq r)].$$
We now use essentially the same proof as in \cite[Lemma 6]{bianeyor}. 
Let $\eps>0$. Since $\tau_{tx}\wedge\eps$ is a stopping time for the
natural filtration of $X$, 
\begin{eqnarray*}
\lefteqn{N^{(1)}[G(t^{1/(1-\alpha)}X_{\tau_{tx}\wedge\eps},0\leq x\leq K)]}\\
&=&E^{\uparrow}[h(X_{
\eps})G(t^{1/(1-\alpha)}X_{\tau_{tx}\wedge\eps},0\leq x\leq K)]\\
&=&
E^{\uparrow}[E^{\uparrow}[h(X_{\eps})|X_{\tau_{tK}\wedge\eps}]G(t^{1/
(1-\alpha)}X_{\tau_{tx}\wedge\eps},0\leq x\leq K)].
\end{eqnarray*}
Since $\tau_{tK}\to 0$ a.s.\ as $t\downarrow 0$, we obtain the same limit if we
remove the $\eps$ in the left-hand side, hence giving 
$\lim N^{(1)}[G(t^{1/(1-\alpha)}L_1^{tx},0\leq x\leq K)]$ by Theorem
\ref{jeulinth}. Using the backwards martingale convergence theorem we
obtain that the conditional expectation on the right-hand side converges to 
$E^{\uparrow}[h(X_{\eps})]=1$. So 
$$\lim_{t\downarrow 0}N^{(1)}[G(t^{1/(1-\alpha)}L_1^{tx},0\leq x\leq
K)]= \lim_{t\downarrow0} E^{\uparrow}[G(t^{1/(1-\alpha)}X_{\tau_{tx}},0\leq
x\leq K)]$$
and the last expression is constant, equal to
$E^{\uparrow}[G(X_{\tau_x},0\leq x\leq K)]$ by scaling, hence the result by
Lamperti's transform. The independence with the initial process is a
refinement of the argument above, using the Markov property at the time
$\tau_{tK}\wedge\eps$.
\cq

\noindent{\bf One final comment. }
It may look quite strange in the proofs above that the {\em a priori} 
ill-defined time $\tau_t$ under the laws $P^{\uparrow}$ or $N$ 
somehow has to be non-degenerate by the proofs we used, even though no 
argument on the path behavior near $0$ has been given for these laws. 
As a matter of fact, things are maybe clearer when considering also the 
inverse Lamperti transform. As above, for some process $Z$ that is 
strictly positive on a set of the form $(0,K)$, $K>0$, we let 
$$C_s=\inf\left\{u\geq 0:\int_0^u \d v\, Z_v>s\right\}.$$
Define the process $X$ by $X_s=Z_{C_s}$. Then we claim that the map
$s \mapsto 1/X_s$ is integrable on a neighborhood of $0$ and that 
$X_{\tau_t}=Z_t$. Indeed, by a change of variables $w=C_v$, one has:
$$\int_0^u \frac{\d v}{X_v}=\int_0^u\frac{\d v}{Z_{C_v}}=\int_0^{C_u}
\frac{Z_w \d w}{Z_w}=C_u<\infty ,$$
as long as $u<C^{-1}(\infty)=\inf\{s:X_s=0\}$, which is 
strictly positive by the hypothesis made on $Z$. This kind of 
arguments also shows that 
as soon as we have one side of Lamperti's theorem, i.e.\ $X_s=Z_{C_s}$ or 
$Z_t=X_{\tau_t}$, with non-degenerate $C$ or $\tau$, 
then the other side is true. In particular, Theorem 
\ref{jeulinth} and Lemma \ref{amauryth} could be restated with 
the inverse statement giving the L\'evy process by time-changing the 
CSBP with $C$. 

\noindent{\bf Acknowledgement. }
I warmly thank Jean Bertoin for suggesting this research and for his many
comments throughout the elaboration of the paper. Thanks also to an 
anonymous referee for useful comments.

\end{document}